\documentclass[10pt]{article}
\usepackage{hyperref}
\usepackage{amsfonts}
\usepackage{mathrsfs}
\usepackage{array}
\usepackage{amsmath,amsthm,amsfonts,amssymb,epsfig}
\usepackage{soul}
\usepackage{color}
\usepackage{stmaryrd}
\usepackage{hyperref}
\soulregister\ref7
\soulregister\bf7

\newtheorem{thm}{Theorem}[section]

\newtheorem{lem}[thm]{Lemma}

\theoremstyle{definition}

\theoremstyle{remark}

\numberwithin{equation}{section}
\usepackage[framemethod=tikz]{mdframed}
\usepackage{xcolor}
\newmdenv[
linecolor=yellow, 
backgroundcolor=yellow!100, 
innerlinewidth=0pt, 
outerlinewidth=0pt, 
leftmargin=0pt, 
rightmargin=0pt
]{highlighted}

\newcommand\be{\begin{equation}}
\newcommand\ee{\end{equation}}
\newcommand\bes{\begin{eqnarray}}
\newcommand\ees{\end{eqnarray}}
\newcommand\bess{\begin{eqnarray*}}
\newcommand\eess{\end{eqnarray*}}

\def\theequation{\arabic{section}.\arabic{equation}}
\begin{document}
	\title{Conditional stability for an inverse problem of a fully discrete stochastic hyperbolic equation}
	
	\author{ Bin Wu$^{1}$,\quad Xu Zhu$^1$, \quad Wenwen Zhou$^1$,\quad Zewen Wang$^2$\footnote{Corresponding author. email: zwwang6@163.com; wangzewen@gzmtu.edu.cn }
		\\
		$^1$School of Mathematics and Statistics\\ Nanjing University of
		Information Science
		and Technology\\
		Nanjing 210044, China\\
        $^2$School of Arts and Sciences\\
        Guangzhou Maritime University\\
        Guangzhou 510725, China
	}

	\maketitle

\begin{abstract}
In this paper, we investigate a discrete inverse problem for a fully discrete approximation of one-dimensional stochastic hyperbolic equation. This inverse problem aims to determine three unknowns, namely the initial displacement, the initial velocity  and the random source term, with discrete spatial derivative data at the left endpoint together with measurements of the solution at the final time. We first prove a new Carleman estimate
for the fully discrete stochastic hyperbolic equation. Based on this Carleman estimate, we then establish stability for this discrete inverse problem under a priori information. Owing to the discrete setting, an extra term depending on the mesh size arises on the right-hand side of the stability estimate.

\vskip 0.3cm


\end{abstract}

	\section{Introduction}
	\setcounter{equation}{0}

In recent decades, inverse problems for stochastic partial differential equations have developed into an important branch of applied mathematics, drawing considerable attention from the scientific community due to their wide range of applications [\ref{Bao-SIAMNA-2020}, \ref{Cav 2002}, \ref{Li-Xu-SIAMUQ-2022}, \ref{Lu2023}]. A pioneering contribution was made in [\ref{Lü 2015}], where global uniqueness of an inverse problem for recovering the initial states and a random source in a stochastic hyperbolic equation was established in the continuous setting by means of Carleman estimates. In this paper, we focus on discrete Carleman estimates and their application to the corresponding discrete inverse problem for fully discrete approximations of one-dimensional stochastic hyperbolic equations.

Let $(\Omega, \mathcal{F}, \{\mathcal{F}_t\}_{t\geq0}, \mathbb P)$ be a complete filtered probability space on which a one-dimensional standard Brownian motion $\{B(t)\}_{t\geq0}$ is defined such that $\{\mathcal F_t\}_{t\geq 0}$ is the natural filtration generated by $B(\cdot)$, augmented by all the $\mathbb P$-null sets in $\mathcal{F}$.  Let $T>0$, $I=(0,1)$ and $Q=I\times(0,T)$.  The continuous stochastic hyperbolic equation can be written as
\begin{align}\label{1.1}
		\left\{
		\begin{array}{ll}
			{\rm d}y_{t}-y_{xx}{\rm d}t=(ay+by_{x}+ c y_t){\rm d}t+(dy+g){\rm d} B(t), &(x,t)\in Q,\\
			y(0,t)=y(1,t)=0, 	&t\in (0,T), \\
			y(x,0)=y_0,\quad y_{t}(x,0)=y_1,		    & x\in I
		\end{array}
		\right.
	\end{align}
with suitable functions $a$, $b$, $c$ and $d$. Physically, the random source $g$ represents the intensity of a random force. Additionally, $y_0$ satisfies the following compatibility condition:
\begin{align}
y_0(0)=y_0(1)=0.
\end{align}

In [\ref{Lü 2015}], the authors investigated the following inverse problem and established its global uniqueness by means of  a Carleman estimate for continuous stochastic hyperbolic equations.

\vspace{2mm}
	
	\noindent{\bf Continuous inverse problem}.\ Determine three unknowns, namely the random source $g$, the initial displacement $y_0$ and the initial velocity $y_1$, from the boundary flux $y_x |_{{x=0}}$ and  the final time data $y(T)$.

\vspace{2mm}

The main objective of this paper is to derive a discrete Carleman estimate for the finite-difference full discretization of the stochastic hyperbolic equation (\ref{1.1}). We then apply
this discrete Carleman estimate to investigate the corresponding discrete version of continuous inverse problem for the fully discrete stochastic hyperbolic equation. To do this,  let us consider $M, N\in \mathbb{N}^*$, the space- and time-discretization parameters $\Delta x=\frac{1}{M+1}$ and $\Delta t=\frac{T}{N}$. Furthermore, we introduce
an equidistant spatial mesh of $(0, 1)$, $0=x_0<x_1<\cdots<x_{M+1}=1$, and an equidistant time mesh of $(0, T)$, $0=t^0<t^1<\cdots<t^N=T$. Let $x_{j} = jh, t^{n}=n\Delta t$, $B^n=B(t^n)$ and 
$$y_{0,j}=y_0(x_j),\quad y_{1,j}=y_1(x_j)\quad {\rm and}\quad   y^n_j=y(x_j,t^n),\quad  0\leq j\leq M+1, \ 0\leq n\leq N,$$
analogous definitions are for $a^n_j, b^n_j, c^n_j, d^n_j, g^n_j$ and the other functions
in the sequel. Then, the explicit finite difference approximation of (1.1) can be written as follows
	\begin{equation}\label{1.2}
		\left\{
		\begin{array}{l}
			\frac{y^{n+1}_j-2{y^{n}_j}+{ y^{n-1}_j}}{\Delta t}-\frac{y^{n}_{j+1}-2y^n_j+y^{n}_{j-1}}{(\Delta x)^2}\Delta t=\left(a^n_j y^n_j+b^{n}_j\frac{y_{j+1}^{n}-y_{j-1}^{n}}{2\Delta x}+c^{n}_j\frac{{ y_j^{n}-y_{j}^{n-1}}}{\Delta t}\right)\Delta t\\
			\hspace{3.2cm}+(d^n_jy^n_j+g^n_j)(B^{n+1}-B^n), \quad  1\leq j\leq M,\ 1\leq n\leq N,\\
			y^{n}_0=y^{n}_{M+1}=0, \hspace{5.0cm} \quad 1\leq n\leq N,\\
			y^{0}_j=y_{0,j}, \quad y^1_{j}=y_j^{0}+y_{1,j}\Delta t ,  \hspace{2.7cm} \quad  0\leq j\leq M+1.
		\end{array}
		\right.
	\end{equation}

To facilitate the formulation of our inverse problem, we introduce some discrete notations. We use the notation 
$\llbracket a, b \rrbracket = [a,b] \cap \mathbb{N}$ 
for any real numbers $a<b$. Let us also define the uniform meshes in space and time as follows:
\begin{align*}
\mathcal K_{\Delta x}=\{x_j:=j\Delta x;\  j\in \llbracket 0, M+1 \rrbracket\}\quad {\rm and}\quad \mathcal K_{\Delta t}=\{t^n:=n\Delta t;\  n\in \llbracket 0, N \rrbracket\}.
\end{align*}
For $x\in \mathcal K_{\Delta x}$, we define the translate operator ${\bf s}_{\pm}(x)=x\pm {\Delta x}/{2}$.
Then for any set of points $\mathcal{W}\subset\mathcal{K}_{\Delta x}$, we define the dual meshes as follows:
	\begin{align*}
	    \mathcal{W}^\prime={\bf s}_{-}(\mathcal{W})\cap{\bf s}_{+}(\mathcal{W}),\quad \mathcal{W}^*={\bf s}_{-}(\mathcal{W})\cup{\bf s}_{+}(\mathcal{W}).
	\end{align*}
    The dual of a dual set will be denoted as $\overline{\mathcal{W}}:=(\mathcal{W}^{*})^{*}$ and $\widetilde{\mathcal{W}}:=(\mathcal{W}')'$.
Then, any subset $\mathcal{W} \subset \mathcal{K}_{\Delta x}$ that verifies $\widetilde{\overline{\mathcal{W}}}=\mathcal{W}$ will be called regular mesh. Furthermore, we define its boundary by $\partial\mathcal{W}:=\overline{\mathcal{W}} \backslash \mathcal{W}$. We denote the interior of the meshes in space by $\mathcal M:= \widetilde {\mathcal K}_{\Delta x}$. We also introduce some notations to define the discretization of the time variable.
We define $\mathcal{N}:=\{t^{n};\ n \in \llbracket 1, N \rrbracket \}$ as the primal mesh in time and $\mathcal{N}^{*}:=\{t^{n+\frac{1}{2}};\ n \in \llbracket 0, N-1 \rrbracket \}$ as the dual mesh in time.
We write $\partial \mathcal N=\{0, T\}$, $\overline{\mathcal{N}}=\mathcal{N} \cup \{0\}$. 

We denote by $\mathbb C(\mathcal M\times \mathcal N)$ the set of functions defined on $\mathcal M\times \mathcal N$. For $u\in \mathbb C(\mathcal M\times \mathcal N)$ we introduce the spatial and temporal translation operators 
\begin{align*}
{\bf s}_{\pm}u (x,t)=u(x\pm\frac{\Delta x}{2},t),\quad {\bf t}^{\pm}u (x,t)=u(x,t\pm \frac{\Delta t}{2})
\end{align*}  
and further define the average and  difference operators as
	\begin{align*}
		&{\bf A}_x u(x,t)=\frac{{\bf s}_+u(x,t)+{\bf s}_-u(x,t)}{2},\quad {\bf D}_x u(x,t)=\frac{{\bf s}_+u(x,t)-{\bf s}_-u(x,t)}{\Delta x},\\
&{\bf A}_t u(x,t)=\frac{{\bf t}^+u(x,t)+{\bf t}^-u(x,t)}{2},\quad {\bf D}_t u(x,t)=\frac{{\bf t}^+u(x,t)-{\bf t}^-u(x,t)}{\Delta t}.
	\end{align*}
Also, we define time-increment operator as
\begin{align*}
{\bf d}_t u(x,t)=\Delta t\, {\bf D}_t u(x,t)={\bf t}^+u(x,t)-{\bf t}^-u(x,t).
\end{align*}

	\noindent{\bf Discrete inverse problem.} \ Determine three unknowns $y_{0}, y_{1}$ and $g$ in the following fully discrete stochastic hyperbolic equation:
	\begin{align}\label{1.3}
		\left\{
		\begin{array}{ll}
			{\bf d}_t z-{\bf D}_x^2y\Delta t=\big(a y+b{\bf A}_x{\bf D}_{x}y+c {\bf t}^{-}(z)\big)\Delta t\\
\hspace{3.9cm}+(dy+g) {\bf t}^+({\bf d}_tB), &(x,t)\in \mathcal M\times \mathcal N,\\
z={\bf D}_t y, &(x,t)\in \mathcal M\times \mathcal N,\\
			y(0,t)=y(1,t)=0, 	& t\in  \mathcal N, \\
			y(x,0)=y_{0},\quad {\bf t}^+(z) (x,0)=y_{1},		   & x\in {\overline {\mathcal M}},
		\end{array}
		\right.
	\end{align}
	by the observation data
	$${\bf tr}({\bf D}_{x}y) \big |_{x=0} \quad{\rm and}\quad y \big |_{t=T},$$
where the trace operator ${\bf tr}$ will be explicitly defined in (\ref{2.2}).

\vspace{2mm}

{The well-posedness of the fully discrete system (\ref{1.3}), along with its discrete energy estimate, will be established in Appendix A.}

\vspace{2mm}

{\noindent\bf Remark 1.1.}\ Unlike the discretizations employed in [\ref{Carreno-ACM-2023}, \ref{Lecaros-JDE-2023}] for studying discrete Carleman estimates of deterministic partial differential equations,  we adopt  an explicit finite difference scheme rather than an implicit one for the  discretization of $y$. Meanwhile, for the Brownian motion $B(t)$, we employ a forward time difference in the temporal discretization. This is because we require the independence between $\mathcal B(y,z)$ and ${\bf t}^+({\bf D}_t B)$ to  guarantee the discrete calculus
\begin{align*}
\mathbb E\iint\limits_{\mathcal M\times\mathcal N}\mathcal B(y, z){\bf t}^+({\bf D}_tB)=\ &\mathbb E\sum_{j=1}^{M}\sum_{n=1}^{N} \mathcal B(y^n_j, z^n_j) (B^{n+1}-B^n)\Delta x\\
=\ &\sum_{j=1}^{M}\sum_{n=1}^{N}\mathbb E(\mathcal B(y^n_j, z^n_j))\mathbb E(B^{n+1}-B^n)\Delta x=0,
\end{align*}
which is essential for establishing  Carleman estimate for stochastic partial differential equations, { details see (\ref{2-A.3}) in Appendix A}. Here, $\mathcal B $ denotes a composition of difference operators involving $y$ and $z$, arising from the discretization of ${\rm d}z-y_{xx}{\rm d}t$.

\vspace{2mm}

{	
It is well known that Carleman estimate is a fundamental tool for studying inverse problems. The ground breaking idea of applications of Carleman estimates to inverse problems was first presented in [\ref{Bukhgeim-SM-1981}], where the uniqueness for an inverse problem was established for deterministic partial differential equations by means of Carleman estimates. Since then, it has become a powerful and widely used methodology for applying Carleman estimates to inverse problems for various types of partial differential equations [\ref{Bellassoued-Springer-2017}, \ref{Imanuvilov-IP-2021}, \ref{Yamamoto-SP-2025}], including discrete inverse problems [\ref{Baudouin-SICON-2013}, \ref{Lecaros-IP-2025}]. Later, studies on both stochastic and discrete inverse problems were also developed upon this foundation. Moreover, Carleman estimates can also be used in the construction of convexification numerical methods for a broad class of coefficient inverse problems, which are of great importance in real applications due to their global convergence property [\ref{Klibanov-DG-2021}, \ref{Klibanov-IPI-2021}].
}

{As for Carleman estimates for stochastic inverse problems, we refer to} [\ref{Dou-IP-2024}, \ref{Lu-IP-2013}, \ref{Lu2023}, \ref{Wu 2022}, \ref{Zha-SIMA-2008}]. For the inverse problem for stochastic hyperbolic equations, Lü [\ref{Lü 2015}] introduced a special weight function to place the unknown random source term on the left-hand side of the Carleman estimate, thereby establishing uniqueness in the determination of three unknowns. This approach was subsequently extended to the stochastic dynamic Euler–Bernoulli beam equation [\ref{Yuan-JMAA-2017}] and the damped stochastic plate equation [\ref{Yu-MMAS-2023}]. We also refer to [\ref{Lu-IP-2012}] for the inverse source problem and [\ref{Yuan-IP-2017}] for the initial data determination problem, both concerning stochastic parabolic equations. Regarding numerical methods for stochastic inverse problems, we mention [\ref{Dou-IP-2024}] for the ill-posed Cauchy problem of stochastic parabolic equation, [\ref{Lu-Wang-2025}] for the inverse source problem of semilinear stochastic hyperbolic equation, and [\ref{Bao2012}] for recovering a random source in quantifying the elastic modulus of nanomaterials. All the above stochastic inverse problems have been studied within a continuous framework.

Let us briefly review relevant studies on discrete Carleman estimates and their applications. Extensive literature exists for deterministic partial differential equations: see [\ref{Bau2015}, \ref{Boyer-JMPA-2010}, \ref{Boyer-POINCARE-AN-2014}, \ref{Cerpa-JMPA-2022}, \ref{Lecaros-ESIAM-2021}] for spatial semi-discrete Carleman estimates, and [\ref{Boyer-ESAIM:COCV-2005}, \ref{Hernández2021}, \ref{Zhang 2022}] for time semi-discrete estimates. In the fully discrete setting, in [\ref{Carreno-ACM-2023}] the authors established a Carleman estimate for implicit finite-difference approximations of one-dimensional parabolic operators, which was then applied to derive relaxed controllability results for fully discrete linear and semilinear parabolic equations. This line of research was subsequently extended to fully discrete parabolic operators with dynamic boundary conditions [\ref{Lecaros-JDE-2023}] and to fully discrete Korteweg–de Vries equations [\ref{Kumar-2024}]. However, discrete Carleman estimates for stochastic partial differential equations remain relatively scarce. Recently, spatial semi-discrete Carleman estimates have been derived for stochastic parabolic equations [\ref{Wu 2024}, \ref{Zhao-SICON-2025}] and successfully applied to controllability and inverse problems. Moreover, a new global Carleman estimate for spatial semi-discrete backward stochastic fourth-order parabolic operators has been used to establish relaxed null controllability for a fourth-order semi-discrete parabolic equation [\ref{Wang-arXiv-2024}].

Relatively few works have focused on discrete inverse problems. For deterministic equations, in [\ref{Baudouin-SICON-2013}] the authors established a stability result for the discrete inverse problem of determining the potential in spatial semi-discrete hyperbolic equations. In [\ref{Klibanov-IP-2024}], ill-posed problems and coefficient inverse problems for deterministic parabolic equations were investigated within a time finite-difference framework. Stability analysis and globally convergent numerical methods for such discrete inverse problems have also been developed. As for discrete stochastic inverse problems, we refer to [\ref{Wu 2024}] for results on spatial semi-discrete stochastic parabolic equations. { These studies on discrete inverse problems have, to some extent, drawn upon the idea introduced in [\ref{Bukhgeim-SM-1981}].}

However, to the best of our knowledge, no existing work has addressed discrete Carleman estimates in conjunction with discrete inverse problems for fully discrete stochastic hyperbolic equations, not even for their fully discrete deterministic counterparts. In this paper, we first establish a Carleman estimate for fully discrete stochastic hyperbolic equations. As an application, this estimate is then employed to study the stability of the associated discrete inverse problem.  Compared with deterministic discrete Carleman estimates, additional terms arising from stochastic effects must also be addressed. { Combining the discrete energy estimate with this Carleman estimate}, we are able to obtain a stability result for the discrete inverse problem under a prior condition on $(g, y_0, y_1)$, which is particularly valuable for error analysis in numerical implementations.

The rest of this paper is organized as follows. Section 2 introduces fundamental concepts and formulas of discrete calculus. In Section 3, we establish a new Carleman estimate for the fully discrete stochastic hyperbolic equation. Section 4 is devoted to deriving  conditional stability result for the corresponding discrete inverse problem. In the last section, we provide some further discussion.

\section{Discrete settings}
\setcounter{equation}{0}

In this section, we introduce fundamental concepts and formulas of discrete calculus, including discrete function spaces, some useful discrete identities, integration by parts for discrete operators and discrete calculus results on the Carleman weight function, which will be utilized in the proofs of the main results. 

\subsection{Discrete function spaces}	

	For $\mathcal W\subset \mathcal K_{\Delta x}$ and $u, v\in \mathbb{C}(\overline {\mathcal W})$, we define the discrete integral of $u$ on $\mathcal{W}$ and $\partial \mathcal W$

\begin{align*}\displaybreak[0]\int\limits_{\mathcal{W}}u
=\sum_{x \in \mathcal{W}}u(x)\Delta x,\quad \int\limits_{\partial\mathcal{W}}u=\sum\limits_{x \in \partial\mathcal{W}}u(x)
\end{align*}
and the following $L^{2}$-inner product
$$(u,v)_{L^{2}(\mathcal{W})}
=\int\limits_{\mathcal{W}}u(x)v(x)
=\sum_{x \in \mathcal{W}}u(x)v(x)\Delta x.$$
 The associated norm is denoted by $\|u\|_{L^2(\mathcal W)}$, i.e.
\begin{align*}
\|u\|_{L^2(\mathcal W)}=\left(\sum\limits_{x \in \mathcal{W}}|u(x)|^2\Delta x\right)^{\frac{1}{2}}.
\end{align*}
 Analogous definitions are for $\|u\|_{L^2(\partial\mathcal W)}$ and $\|u\|_{L^{\infty}(\mathcal{W})}$ in terms of
$$\|u\|_{L^2(\partial\mathcal W)}=\left(\sum\limits_{x\in \partial\mathcal{W}}|u(x)|^2\right)^{\frac{1}{2}},\quad \|u\|_{L^{\infty}(\mathcal{W})}=\max \limits_{x \in \mathcal{W}}|u(x)|.$$
Also, we introduce
\begin{align*}
\|u\|_{H^1(\mathcal W)}=\bigg(\ \int\limits_{\mathcal W^*} |\mathbf D_x u|^2+\int\limits_{\mathcal W}|u|^2\bigg)^{\frac{1}{2}}.
\end{align*}

Similarly, for $v\in \mathbb{C}(\overline {\mathcal N})$, we define the discrete integral of $v$ with respect to time on $\mathcal N$ and $\partial\mathcal N$
$$\int\limits_{\mathcal{N}} v
= \sum_{t \in \mathcal{N}}v(t)\Delta t,\quad \int\limits_{\partial\mathcal{N}} v
=\sum_{t \in \partial\mathcal{N}}v(t).$$
Further, we introduce
\begin{align*}
\|v\|_{L^2(\mathcal N)}=\bigg(\ \int\limits_{\mathcal{N}} |v|^2\bigg)^{\frac{1}{2}},\quad \|{\bf D}_t v\|_{L^2(\mathcal N)}=\bigg(\ \int\limits_{\mathcal{N}} {\bf t}^-\left(|{\bf D}_t v|^2\right)\bigg)^{\frac{1}{2}}.
\end{align*}

To introduce the boundary conditions, for $x \in \partial\mathcal{W}$ we define the outward normal at $x$ as
\begin{equation}\label{}
n_x(x):=
	\left\{
\begin{aligned}
&1,
	&&{\bf s}_{-}(x) \in \mathcal{W}^{*} \; \text{and} \; {\bf s}_{+}(x) \notin \mathcal{W}^{*},\\
-&1,
	&&{\bf s}_{-}(x) \notin \mathcal{W}^{*} \; \text{and} \; {\bf s}_{+}(x) \in \mathcal{W}^{*},\\
&0,
	&&\text{otherwise}\\
\end{aligned}
\right.
\end{equation}
and the trace operator for $u \in \mathbb{C}(\mathcal{W}^{*})$ as
\begin{equation}\label{2.2}
{\bf tr}(u)(x):=
	\left\{
\begin{aligned}
&{\bf s}_{-}u(x),
	&&n_x(x)=1,\\
&{\bf s}_{+}u(x),
	&&n_x(x)=-1,\\
&0,
	&&n_x(x)=0.\\
\end{aligned}
\right.
\end{equation}
Similarly,  we define the outward normal for $t\in \partial \mathcal N$ as
\begin{equation}\label{1.6}
n_t(t):=
	\left\{
\begin{aligned}
&1,
	&&t=T,\\
-&1,
	&&t=0.
\end{aligned}
\right.
\end{equation}

Combining space and time discrete integrals, for $u\in \mathbb{C}(\mathcal K_{\Delta x}\times\mathcal K_{\Delta t})$ we define 
\begin{align*}
&\iint\limits_{\mathcal M\times\mathcal N}u=\sum_{x \in \mathcal{M}}\sum_{t\in \mathcal{N}} u(x,t)\Delta x\Delta t,\quad \iint\limits_{\partial\mathcal M\times\mathcal N}u=\sum_{x \in \partial\mathcal{M}}\sum_{t \in \mathcal{N}} u(x,t)\Delta t,\\
&\iint\limits_{\mathcal M\times\partial\mathcal N}u=\sum_{t \in \partial\mathcal{N}}\sum_{x \in \mathcal{M}}u(x,t)\Delta x,\quad \iint\limits_{\partial\mathcal M\times\partial\mathcal N}u=\sum_{x \in \partial\mathcal{M}}\sum_{t \in \partial\mathcal{N}} u(x,t),
\end{align*}
and the $L^2$-inner product in $\mathcal M\times \mathcal N$
\begin{align*}(u,v)_{L^{2}(\mathcal{M}\times\mathcal N)}
=\iint\limits_{\mathcal{M}\times\mathcal N}uv
=\sum_{x \in \mathcal{M}}\sum_{t\in \mathcal N}u(x,t)v(x,t)\Delta x\Delta t.
\end{align*}
The associated norm is denoted by $\|u\|_{L^2(\mathcal M\times\mathcal N)}$, i.e.
\begin{align*}
\|u\|_{L^2(\mathcal M\times \mathcal N)}=\left(\sum\limits_{t \in \mathcal{N}}\sum\limits_{x \in \mathcal{M}}|u(x, t)|^2\Delta x\Delta t\right)^{\frac{1}{2}}.
\end{align*}

Additionally, we introduce discrete stochastic function spaces for stochastic analysis on discrete meshes. For a discrete Banach space $\mathcal Y$, we denote by $L_\mathcal F^2(\Omega;\mathcal Y)$ the space of all progressively measurable stochastic process $\zeta$ such that $\mathbb E(\|\zeta\|^2_{\mathcal Y})<\infty$, e.g.
\begin{align*}
\|\zeta\|_{L^2_{\mathcal F}(\Omega;L^2(\mathcal M\times\mathcal N))}=\mathbb E \left(\sum\limits_{t \in \mathcal{N}}\sum\limits_{x \in \mathcal{M}}|\zeta(x, t)|^2\Delta x\Delta t\right)^{\frac{1}{2}}.
\end{align*}
Analogous definition is for $L_\mathcal F^\infty(\Omega; \mathcal Y)$.

\subsection{Discrete calculus formulas}

For the spatial difference operators ${\bf A}_x$ and ${\bf D}_x$, we provide several preliminary identities together with discrete integration by parts formulas that will be frequently used in the subsequent analysis. The corresponding proofs can be found in [\ref{Boyer-JMPA-2010}].

\begin{lem}  Let $ u, v\in \mathbb C(\overline {\mathcal W})$ with $\mathcal W\subset \mathcal K_{\Delta x}$. Then for the average operator ${\bf A}_x$ and difference operator ${\bf D}_x$, we have the following identities:
\begin{align}
			\label{}& {\bf A}_x(uv)={\bf A}_x u{\bf A}_x v+\frac{1}{4}(\Delta x)^2{\bf D}_xu{\bf D}_xv,\\
			\label{}&\quad\ \  {\bf D}_x(uv)={\bf D}_x u{\bf A}_x v+{\bf A}_x u{\bf D}_x v
\end{align}
and further
\begin{align}\label{2.6}
	u={\bf A}_x^2 u-\frac{1}{4}(\Delta x)^2{\bf D}_x^2 u.
\end{align}
	\end{lem}
	
\begin{lem} Let $ u\in \mathbb C(\overline{\mathcal W}) $ and $v\in \mathbb C(\mathcal W^*)$ with $\mathcal W\subset \mathcal K_{\Delta x}$. Then we have the following identities:
\begin{align} 
			\label{2.7}&\int\limits_{\mathcal W}u {\bf A}_xv=\int\limits_{\mathcal W^*}{\bf A}_{x}u v-\frac{1}{2}\Delta x\int\limits_{\partial\mathcal W}u {\bf tr}(v),\\
			\label{2.8}&\int\limits_{\mathcal W}u {\bf D}_xv=-\int\limits_{\mathcal W^*}{\bf D}_{x}uv+\int\limits_{\partial \mathcal W}u{\bf tr}(v)n_x.		
\end{align}
\end{lem}

As for time discrete calculus, we have the following formulas [\ref{Boyer-ESAIM:COCV-2005}] or [\ref{Lecaros-JDE-2023}].
 
\begin{lem}
Let $f, g\in \mathbb{C}(\overline {\mathcal N})$. Then for the difference operator ${\bf D}_t$, we have the following identities
\begin{align}
\label{2.9} &\qquad\ \ {\bf D}_{t}(fg)={\bf D}_{t}f{\bf t}^{-}(g)+{\bf t}^{+}(f){\bf D}_{t}g={\bf D}_{t}f{\bf t}^{+}(g)+{\bf t}^{-}(f){\bf D}_{t}g,\\
\label{2.10} &2{\bf t}^{+}(f){\bf D}_{t}f={\bf D}_{t}(f^{2})+\Delta t ({\bf D}_{t}f)^{2}, \quad 2{\bf t}^{-}(f){\bf D}_{t}f={\bf D}_{t}(f^{2})-\Delta t ({\bf D}_{t}f)^{2}.
\end{align}
\end{lem}


%

\begin{lem}
Let $f\in \mathbb{C}(\overline{\mathcal{N}})$ and $g \in \mathbb{C}(\overline{\mathcal{N}^{*}})$.
Then the following identities hold
\begin{align}
\label{2.11}&\int\limits_{\mathcal{N}} f{\bf t}^{-}(g)
=\int\limits_{\mathcal{N}^*}{\bf t}^{+}(f)g,\quad \int\limits_{\mathcal{N}} f{\bf t}^{+}(g)
=\int\limits_{\mathcal{N}^*}{\bf t}^{-}(f)g+\Delta t\int\limits_{\partial \mathcal N} f{\bf t}^+(g)n_t,\\
\label{2.12}&\qquad\qquad\quad\quad\ \ \int\limits_{\mathcal{N}} f{\bf D}_{t}g
=-\int\limits_{\mathcal{N}^*} g{\bf D}_{t}f+\int\limits_{\partial\mathcal{N}}f{\bf t}^{+}(g)n_{t},
\end{align}
where $\partial\mathcal{N}=\{0,T\}$.
\end{lem}

Furthermore, we have the following useful identity, which can be found in [\ref{Lecaros-JDE-2023}]
\begin{align}
\label{2.13}&\int\limits_{\mathcal{N}}{\bf t}^{-}(f){\bf D}_{t}g
=-\int\limits_{\mathcal{N}}{\bf D}_{t}f{\bf t}^{+}(g)+\int\limits_{\partial\mathcal{N}} {\bf t}^{+}(fg)n_{t}
\end{align}
for $f,g \in \mathbb{C}(\overline{\mathcal{N}^{*}})$, and
\begin{align}
\label{1-2.14}
&\int\limits_{\mathcal{N}^*}{\bf t}^{+}(f){\bf D}_{t}g
=-\int\limits_{\mathcal{N}^*}{\bf t}^{-}(g){\bf D}_{t}f+\int\limits_{\partial \mathcal{N}} fgn_{t}
\end{align}
for $f,g \in \mathbb{C}(\overline{\mathcal{N}})$.

%

\subsection{Weight function and its discrete calculus}

We first introduce a weight function that will be used in our discrete Carleman estimate. Let $s$ and $\lambda$ 
be two large parameters. For some $x^*>1$ and a small positive parameter $\beta$,  we define the regular weight function by 
\begin{equation}\label{2.15}
	r=e^{l(x,t)},\quad l(x,t)=s\varphi (x,t),\quad \varphi(x,t)=e^{\lambda \phi(x,t)},\quad  (x,t)\in Q
\end{equation}
with
\begin{equation}\label{}
	\phi(x,t)=|x-x^*|^2-\beta {|t-T|^{2}}+M,\quad (x,t)\in Q,
\end{equation} 
where $M$ is a large positive constant so that $\phi(x,t)>0$ for all $(x,t)\in Q$. We also set $\rho=r^{-1}$.

\vspace{2mm}

     Next, we present some discrete calculus results for the Carleman weight function without 
proofs. We refer to [\ref{Boyer-ESAIM:COCV-2005}, \ref{Boyer-JMPA-2010}] and [\ref{Lecaros-JDE-2023}] for a complete discussion. In contrast to the weight function used for discrete parabolic operator in Carleman estimates,  the corresponding weight for hyperbolic equation possesses sufficiently smoothness in time variable. This property ensures that analogous results hold for the time difference operators ${\bf A}_t$ and ${\bf D}_t$.

\begin{lem} \label{gj2}
Let $\alpha,k,l,m,n \in \mathbb{N}$.
Provided $\max\{{s \Delta x}, s\Delta t\} \le 1$, we have
\begin{align*}
 \ &\ \ \ r{\bf A}_{\xi}^{m}{\bf D}_{\xi}^{n}\partial_{\xi}^{\alpha}\rho= r\partial_{\xi}^{n+\alpha}\rho+s^{n+\alpha}{\mathcal 
O}_\lambda\left((s\Delta \xi)^{2}\right)=s^{n+\alpha}\mathcal O_\lambda(1),\\
\ &{\bf A}_{\xi}^{k}{\bf D}_{\xi}^{l}\partial_{\xi}^{\alpha}(r{\bf A}_{\xi}^{m}{\bf D}_{\xi}^{n}\rho) 
=\partial_{\xi}^{l+\alpha}(r\partial_{\xi}^{n}\rho)+s^{n}\mathcal O_\lambda\left((s\Delta \xi)^{2}\right)=s^{n}{\mathcal O}_\lambda(1),\quad \xi\in\{x,t\}.
\end{align*}
\end{lem}

We also need the following asymptotic expressions of mixed difference operators.

\begin{lem}
Let $m, n\in \mathbb{N}$. Provided $\max\{{s \Delta x}, s\Delta t\} \le 1$ and $\Delta t=\mathcal O(\Delta x)$, we have
\begin{align*}
\ &{\bf A}_t^m{\bf D}_t^n(r{\bf A}_x{\bf D}_x\rho)=\partial_t^n(r\partial_x\rho)+\mathcal O_{\lambda}(s^2\Delta x),\quad {\bf A}_x^m{\bf D}_x^n(r{\bf A}_t{\bf D}_t\rho)=\partial_x^n(r\partial_t\rho)+\mathcal O_{\lambda}(s^2\Delta x),\\
\ &\ \ \ {\bf A}_t^m{\bf D}_t^n(r{\bf D}_x^2\rho)=\partial_t^n(r\partial_x^2\rho)+\mathcal O_{\lambda}(s^3\Delta x),\quad {\bf A}_x^m{\bf D}_x^n(r{\bf D}_t\rho)=\partial_x^n(r\partial_t\rho)+\mathcal O_{\lambda}(s^2\Delta x).
\end{align*}
\end{lem}

{\noindent\bf Proof.}\ We prove only the first result here. The remaining results follow by similar arguments. 

Obviously, we have
\begin{align}
{\bf A}_x{\bf D}_x\rho={\bf A}_x(\partial_x\rho+R_{{\bf D}_x}(\rho))=\partial_x \rho+R_{{\bf A}_x}(\partial_x\rho)+R_{{\bf D}_x}(\rho)+R_{{\bf A}_x{\bf D}_x}(\rho).
\end{align}
Then we obtain from Proposition 4.1 in [\ref{Lecaros-ESIAM-2021}] that
\begin{align}\label{2.19}
&{\bf D}_t^n (r{\bf A}_x{\bf D}_x\rho)\nonumber\\
=\ & \sum_{k=0}^n\binom{n}{k}\partial_t^k\left(R_{{\bf D}_t^{n-k}}\left(r\partial_x\rho+rR_{{\bf A}_x}(\partial_x\rho)+rR_{{\bf D}_x}(\rho)+rR_{{\bf A}_x{\bf D}_x}(\rho)\right)\right),
\end{align}
where
\begin{align*}
R_{{\bf D}_\xi^k}(f)=(\Delta \xi)^2\sum_{l=0}^k\binom{k}{l}(-1)^{l}\left(\frac{k-2l}{2}\right)^{k+2}\int_0^1\frac{(1-\sigma)^{k+1}}{(k+1)!} f^{(k+2)}\Big(\xi+(k-2l)\sigma\frac{\Delta\xi}{2}\Big){\rm d}\sigma,\nonumber\\
R_{{\bf A}_\xi^k}(f)=\frac{1}{2^{k+2}}(\Delta \xi)^2\sum_{l=0}^k\binom{k}{l}(k-2l)^2\int_0^1(1-\sigma) f^{(2)}\Big(\xi+(k-2l)\sigma\frac{\Delta\xi}{2}\Big){\rm d}\sigma,\quad\xi\in \{x,t\}
\end{align*}
and $R_{{\bf A}_\xi^m{\bf D}_\xi^n}=R_{{\bf A}_\xi^m}\circ R_{{\bf D}_\xi^n}$.

Let $f_{\sigma,\tau}(x,t)=f(x+\sigma \frac{\Delta x}{2},t+\tau \frac{\Delta t}{2})$. We easily see that
\begin{align*}
&\partial_t^k\left(R_{{\bf D}_t^{n-k}}(r\partial_x\rho)\right)\nonumber\\
=\ & (\Delta t)^2\sum_{l=0}^{n-k}\binom{n-k}{l}(-1)^{l}\left(\frac{n-k-2l}{2}\right)^{n-k+2}\int_0^1\frac{(1-\tau)^{n-k+1}}{(n-k+1)!} \partial_t^{n+2}(r\partial_x\rho)_{0, (n-k-2l)\tau}{\rm d}\tau,
\end{align*}
together with $\partial_t^{n+2}(r\partial_x\rho)=\partial_t^{n+2}(-s\lambda\varphi\partial_x\phi)=s\mathcal O_{\lambda}(1)$, which implies
\begin{align}\label{2.20}
&\partial_t^k\left(R_{{\bf D}_t^{n-k}}(r\partial_x\rho)\right)=
\left\{\begin{array}{ll}
s\mathcal O_{\lambda}\left((\Delta t)^2\right),& k=0,1,2,\cdots, n-1,\\
\partial_t^n(r\partial_x\rho),& k=n.
\end{array}
\right.
\end{align}
Similarly, we obtain 
\begin{align}\label{2.21}
\left\{\begin{array}{l}
\partial_t^k\left(R_{{\bf D}_t^{n-k}}(rR_{{\bf A}_x}(\partial_x\rho))\right)=s^3\mathcal O_{\lambda}\left((\Delta x\Delta t)^2\right),\\  
\partial_t^k\left(R_{{\bf D}_t^{n-k}}\left(rR_{{\bf D}_x}(\rho)\right)\right)=s^3\mathcal O_{\lambda}\left((\Delta x\Delta t)^2\right),\\
\partial_t^k\left(R_{{\bf D}_t^{n-k}}\left(rR_{{\bf A}_x{\bf D}_x}(\rho)\right)\right)=s^5\mathcal O_{\lambda}\left((\Delta x)^4(\Delta t)^2\right)
\end{array}
\right.
\end{align}
for $k=0, 1,2,\cdots, n-1$  and
\begin{align}\label{2.22}\left\{\begin{array}{l}
\partial_t^k\left(R_{{\bf D}_t^{n-k}}(rR_{{\bf A}_x}(\partial_x\rho))\right)=s^3\mathcal O_{\lambda}\left((\Delta x)^2\right),\\  
\partial_t^k\left(R_{{\bf D}_t^{n-k}}\left(rR_{{\bf D}_x}(\rho)\right)\right)=s^3\mathcal O_{\lambda}\left((\Delta x)^2\right),\\
\partial_t^k\left(R_{{\bf D}_t^{n-k}}\left(rR_{{\bf A}_x{\bf D}_x}(\rho)\right)\right)=s^5\mathcal O_{\lambda}\left((\Delta x)^4\right)
\end{array}
\right.
\end{align}
for $k=n$.
Then substituting (\ref{2.20})-(\ref{2.22}) into (\ref{2.19}) yields
\begin{align*}
{\bf D}_t^n(r{\bf A}_x{\bf D}_x\rho)=\partial_t^n(r\partial_x\rho)+\mathcal O_{\lambda}(s^2\Delta x).
\end{align*}
Therefore, we obtain
\begin{align*}
{\bf A}_t^m{\bf D}_t^n(r{\bf A}_x{\bf D}_x\rho)=\ &\partial_t^n(r\partial_x\rho)+R_{{\bf A}_t^m}(\partial_t^n(r\partial_x\rho))+\mathcal O_{\lambda}(s^2\Delta x)\\
=\ &\partial_t^n(r\partial_x\rho)+\mathcal O_{\lambda}(s(\Delta t)^2)+\mathcal O_{\lambda}(s^2\Delta x)\\
=\ &\partial_t^n(r\partial_x\rho)+\mathcal O_{\lambda}(s^2\Delta x).
\end{align*}
The last three results can follow from similar arguments. \hfill$\Box$

\section{Discrete Carleman estimates}\numberwithin{equation}{section} 
	
This section is devoted to proving Carleman estimate for the fully discrete stochastic hyperbolic equation
\begin{align}\label{3.1}
		\left\{
		\begin{array}{ll}
			{\bf d}_t z-{\bf D}_x^2y\Delta t=f\Delta t+g{\bf t}^+({\bf d}_t B), &(x,t)\in \mathcal M\times \mathcal N,\\
z={\bf D}_t y, &(x,t)\in \mathcal M\times \mathcal N,\\
			y(0,t)=y(L,t)=0, 	&t\in  {\mathcal N}, \\
			y(x,0)=y_{0},\quad {\bf t}^+(z) (x,0)=y_{1},		    &x\in {\overline {\mathcal M}}.
		\end{array}
		\right.
	\end{align}

The proof of our Carleman estimate builds upon the fundamental framework established in the continuous setting [\ref{Lü 2015}, \ref{Zha-SIMA-2008}]. However, discretization in both space and time gives rise to additional terms. In particular, the treatment of discrete computations involving random terms introduces specific technical challenges. To overcome these difficulties, we adapt techniques from existing works on semi-discrete and fully discrete partial differential operators [\ref{Boyer-ESAIM:COCV-2005}, \ref{Carreno-ACM-2023}, \ref{Lecaros-JDE-2023}].


The main result in this section is the following discrete Carleman estimate.
	
	\begin{thm}
		Let $ f, g\in L^2_\mathcal F(\Omega; L^2(\mathcal M\times \mathcal N))$ {satisfy $g(T)\in  L^2_\mathcal F(\Omega; L^2(\mathcal M))$}, $y_0\in L^2_\mathcal F(\Omega; H^1$ $(\mathcal M))$, $y_1\in L^2_\mathcal F(\Omega; L^2(\mathcal M))$, $\beta\in (0,1)$, and let $T$ satisfy
\begin{align}\label{1-3.2}
T>\frac{\sup_{x\in I}|x-x^*|}{\beta}.
\end{align} 
Then for all fixed $\lambda\geq \lambda_0$ with a sufficient large positive constant $\lambda_0$ depending on $x^*, \beta, M,  T$, there exist positive constant $C$ depending on $x^*, \beta, M,  T$, and $C(\lambda)$, $s_{0}>1$, $0<\varepsilon<1$, $0<\widehat {\Delta x}<1$ also depending on $\lambda$  such that 
\begin{align}\label{3.2}
		\ &\mathbb{E} \iint\limits_{\mathcal M\times\mathcal N}s^3\lambda^3\varphi^3 r^2|y|^2+{\mathbb{E} \iint\limits_{\mathcal M\times\mathcal N^*}s\lambda\varphi  r^2|{\bf D}_t y|^2}+\mathbb{E} \iint\limits_{\mathcal M^*\times\mathcal N}s\lambda\varphi r^2 |{\bf D}_xy|^2+ \nonumber\\
		&\mathbb{E} \iint\limits_{\mathcal M\times\mathcal N}s\lambda\varphi{ (T-t)} r^2|g|^2+\mathbb{E} \int\limits_{\mathcal M}s^3\lambda^3\varphi^3r^2|y_0|^2+\mathbb{E} \int\limits_{\mathcal M^*}s\lambda \varphi r^2|{\bf D}_x y_0|^2+\mathbb{E} \int\limits_{\mathcal M}s\lambda \varphi r^2|y_1|^2\nonumber\\
\leq \ & C \mathbb E\iint\limits_{\mathcal M\times\mathcal N}r^2|f|^2+C\mathbb E\int\limits_{\mathcal N}s\lambda\varphi{\bf tr}(|{\bf D}_xy|^2)\bigg|_{x=0}+C(\Delta x)^2\mathbb E\iint\limits_{\mathcal M^*
\times\mathcal N^*}s\lambda^2\varphi|{\bf D}_t{\bf D}_x y|^2+\nonumber\\
\ & C(\lambda)s^3e^{C(\lambda)s}\|y(T)\|_{L^2_\mathcal F(\Omega; H^1(\mathcal M))}+{ C(\Delta t)^{\frac{1}{2 }} \mathbb E\int\limits_{\mathcal M}r^2|g|^2\bigg|_{t=T}}
	\end{align}	
		for all $s\geq s_0$, $0<\Delta x\leq \widehat {\Delta x} $, $0<\Delta t\leq 1$ satisfying the conditions
\begin{align*}
s\Delta x \leq \varepsilon\quad{\rm and}\quad \Delta t=\varepsilon \mathcal O\left((\Delta x)^2\right).
\end{align*}
\end{thm}

\noindent{\bf Remark 3.1.}\ For fully discrete deterministic parabolic equations, the following conditions on $\Delta x$ and $\Delta t$ are imposed to prove discrete Carleman estimates [\ref{Carreno-ACM-2023}, \ref{Lecaros-JDE-2023}]:
\begin{align*}
\frac{s^4\Delta t}{\delta ^4 T^6}\leq \varepsilon\quad {\rm and}\quad \frac{s\Delta x}{\delta T^2}\leq \varepsilon.
\end{align*} 
These conditions imply that $\Delta t=\mathcal O(s^{-4})$, which is more restrictive than our condition $\Delta t=\mathcal O(s^{-2})$. In other words, for discrete hyperbolic equations, the requirement of $\Delta t$ in the proof of Carleman estimate can be relaxed compared to that for discrete parabolic equations. Additionally, the parameter $\varepsilon$ in $\Delta t=\varepsilon \mathcal O\left((\Delta x)^2\right)$ is introduced to handle terms involving random ${\rm d} z$. However, for the fully discrete deterministic hyperbolic equation,   this $\varepsilon$ can be removed from $\Delta t=\varepsilon \mathcal O\left((\Delta x)^2\right)$ when proving Carleman estimate.

\vspace{2mm}
	
	\noindent{\bf Remark 3.2.}\ Similar to spatial semi-discrete deterministic hyperbolic equation [\ref{Baudouin-SICON-2013}], the discrete Carleman estimate (\ref{3.2}) has an additional term of $|{\bf D}_t{\bf D}_{x}y|$ on the right-hand side,  compared with its continuous counterpart in [\ref{Lu-IP-2013}].  This term essentially arises from the spatial discretization.  Although it can not be removed in discrete Carleman estimate, it remains compatible with the Lax argument, details see [\ref{Baudouin-SICON-2013}].  {  The error term associated with $g$ arising from time discretization, namely  $
C(\Delta t)^{\frac{1}{2}}\mathbb E\int_{\mathcal M}r^2|g|^2\big|_{t=T}$,  
can not be absorbed by $\mathbb E\iint_{\mathcal M\times\mathcal N}s \lambda\varphi (T-t) r^2|g|^2$.  
 Consequently, in comparison with the continuous counterpart in [\ref{Lü 2015}], our Carleman estimate entails an additional term $g$ localized at $t=T$ on the right-hand side.}  Nevertheless, in the limit as $\Delta x$ and $\Delta t$ approach zero, under suitable regularity we can eliminate these terms and then obtain the corresponding continuous Carleman estimate.

\vspace{2mm}
	
\noindent{\bf Remark 3.3.}\ In this Carleman estimate, the boundary observation at $x_{0}$ can be replaced by the one at $x_{N+1}$ if we choose $x^*<0$ to satisfy $\partial_x\phi>0$ for all $x\in I$.

\vspace{2mm}

\noindent{\bf Remark 3.4.}\ Under the same assumptions in Theorem 3.1,  for fully discrete deterministic hyperbolic equation (i.e. $g=0$),   we can obtain the following Carleman estimate:
\begin{align*}
		\ &\iint\limits_{\mathcal M\times\mathcal N}s^3\lambda^3\varphi^3 r^2|y|^2+\iint\limits_{\mathcal M\times\mathcal N^*}s\lambda\varphi  r^2|{\bf D}_t y|^2+\iint\limits_{\mathcal M^*\times\mathcal N}s\lambda\varphi r^2 |{\bf D}_xy|^2+ \nonumber\\
		&\int\limits_{\mathcal M}s^3\lambda^3\varphi^3r^2|y_0|^2+\int\limits_{\mathcal M^*}s\lambda \varphi r^2|{\bf D}_x y_0|^2+\int\limits_{\mathcal M}s\lambda \varphi r^2|y_1|^2\nonumber\\
\leq \ & C \iint\limits_{\mathcal M\times\mathcal N}r^2|f|^2+C\int\limits_{\mathcal N}s\lambda\varphi{\bf tr}(|{\bf D}_x y|^2)\bigg|_{x=0}+C(\Delta x)^2\iint\limits_{\mathcal M^*
\times\mathcal N^*}s\lambda^2\varphi|{\bf D}_t{\bf D}_x y|^2,
	\end{align*}
where the terminal condition $y(T)=0$ is imposed, and this condition is achievable by a standard cut-off argument, see [\ref{Bellassoued-Springer-2017}]. This Carleman estimate has potential application to reconstruction of initial conditions of fully discrete hyperbolic equations by additional boundary observation.
	
\vspace{2mm}
	
{\noindent\bf Proof of Theorem 3.1.}\ 	 The whole proof  of the Carleman estimate (\ref{3.2})  is divided into the following four steps.
	 
{\em Step 1.  The change of variable for fully discrete stochastic hyperbolic operator}

	Firstly, we set
\begin{align}\label{3.3}
	\mathcal{L} (y):={\bf d}_t z-{\bf D}_x^2y\Delta t,\quad z={\bf D}_t y.
\end{align}
We introduce the change of variable $y=\rho Y=r^{-1} Y$ with $r$ given by $(\ref{2.15})$. A direct calculation gives
\begin{align}
   \label{3.4} z=\ &{\bf D}_t (\rho Y)={\bf t}^+(\rho) Z+{\bf D}_t\rho{\bf t}^- (Y),\\
    \label{3.5}{\bf d}_t z=\ &\rho{\bf d}_tZ+{\bf t}^+({\bf d}_t\rho){\bf t}^+(Z) +{\bf t}^-({\bf d}_t\rho){\bf t}^-(Z) +{\bf d}_t
({\bf D}_t\rho)Y,\\
    \label{3.6}{\bf D}_x^2y=\ &{\bf D}_x^2(\rho Y)={\bf D}_x^2\rho {\bf A}_x^2Y+2{\bf A}_x{\bf D}_x\rho {\bf A}_x{\bf D}_xY
+{\bf A}_x^2\rho {\bf D}_x^2Y
\end{align}
with $Z={\bf D}_t Y$.	Together with (\ref{2.6}), substituting (\ref{3.4})-(\ref{3.6}) into (\ref{3.3}) yields	
\begin{align}\label{3.7}
	r\mathcal{L}(y)=\ &{\bf d}_tZ+r{\bf t}^+({\bf D}_t\rho){\bf t}^+(Z)\Delta t +r{\bf t}^-({\bf D}_t\rho){\bf t}^-(Z)\Delta t+r{\bf D}^2_t\rho Y\Delta t-\nonumber\\
	&r{\bf D}_x^2\rho Y\Delta t-2r{\bf A}_x{\bf D}_x\rho {\bf A}_x{\bf D}_xY\Delta t-r\left({\bf A}_x^2\rho+\frac{1}{4}(\Delta x)^2{\bf D}_x^2 \rho\right) {\bf D}_x^2Y\Delta t.
\end{align}
On the other hand, we easily see that
\begin{align*}
r{\bf t}^+({\bf D}_t\rho){\bf t}^+(Z)+r{\bf t}^-({\bf D}_t\rho){\bf t}^-(Z)=2r{\bf A}_t{\bf D}_t\rho {\bf t}^-(Z)+r{\bf t}^+({\bf D}_t\rho){\bf d}_tZ.
\end{align*}
Then, (\ref{3.7}) can be rewritten as
\begin{align}\label{3.8}
	r\mathcal{L}(y)={\bf d}_t Z+\mathcal A(Y)\Delta t+\mathcal B(Y,Z)\Delta t-\mathcal R(Y,Z)\Delta t,\quad Z={\bf D}_t Y,
\end{align}
where
\begin{align*}
\mathcal A(Y)=\sum_{i=1}^2 \mathcal A_i(Y),\quad \mathcal B(Y,Z)=\mathcal B_1(Z)+\sum_{i=2}^3 \mathcal B_i(Y)
\end{align*}
with
\begin{align*}
    &\mathcal A_1(Y)=-r\left({\bf A}_x^2\rho +\frac{1}{4}(\Delta x)^2{\bf D}_x^2 \rho\right){\bf D}_x^2 Y,\\
    &\mathcal A_2(Y)=\left(r{\bf D}_t^2 \rho-r{\bf D}_x^2\rho \right)Y,\\
	&\mathcal B_{1}(Z)=2r{\bf A}_t{\bf D}_t\rho {\bf t}^-(Z),\\
    &\mathcal B_2(Y)=-2r {\bf A}_x{\bf D}_x\rho {\bf A}_x{\bf D}_x Y,\\ 
    &\mathcal B_3(Y)=r\varphi^{-1}\left({\bf D}_t\varphi {\bf D}_t\rho-{\bf D}_x\varphi {\bf D}_x\rho \right)Y,\\
\end{align*}
and 
\begin{align*}
\mathcal R(Y,Z) =\ & r\varphi^{-1}\left({\bf D}_t\varphi {\bf D}_t\rho-{\bf D}_x\varphi {\bf D}_x\rho\right)Y-r{\bf t}^+({\bf D}_t\rho){\bf d}_tZ.
\end{align*}
    Furthermore, multiplying both sides of (\ref{3.8}) by $2 \mathcal B(Y, Z)$, it follows that
\begin{align}\label{3.9}
	&2r\mathcal{L}(y)\mathcal B(Y, Z)\nonumber\\
 =\ &2\mathcal B(Y, Z){\bf d}_t Z+2\mathcal A(Y)\mathcal B(Y, Z)\Delta t+2|\mathcal B(Y, Z)|^2\Delta t-2\mathcal R(Y,Z) 
     \mathcal B(Y,Z) \Delta t.
\end{align}
	 Then, integrating  (\ref{3.9}) over $\mathcal M\times \mathcal N$, applying the equation in (\ref{3.1}) and taking mathematical expectation, we find that
\begin{align}\label{3.10}
	&2\mathbb E\iint\limits_{\mathcal M\times \mathcal N} r f \mathcal B(Y,Z)+2\mathbb E\iint\limits_{\mathcal M\times 
     \mathcal N} \mathcal R(Y,Z)\mathcal B(Y,Z)\nonumber\\
 =\ &2\mathbb E\iint\limits_{\mathcal M\times \mathcal N}\mathcal B(Y, Z){\bf D}_t Z+2\mathbb E\iint\limits_{\mathcal M\times \mathcal N}\mathcal A(Y)\mathcal B(Y,Z)+2\mathbb E\iint\limits_{\mathcal M 
     \times \mathcal N} |\mathcal B(Y,Z)|^2,
\end{align}
where we have used that 
\begin{align}\label{1-3.12} 2\mathbb E \iint\limits_{\mathcal M\times \mathcal N}rg\mathcal B(Y,Z){\bf t}^+ ({\bf D}_t B)=0.
\end{align}
	By applying Young's inequality to the two terms on the left-hand side of (\ref{3.10}), we obtain
\begin{align}\label{1-3.12}
	     &\mathbb{E}\iint\limits_{\mathcal M\times\mathcal N}r^2 |f|^2+\mathbb{E}\iint\limits_{\mathcal M\times \mathcal N}|\mathcal R(Y,Z) |^2 \nonumber\\
    \ge\ & \mathbb E\iint\limits_{\mathcal M\times \mathcal N}\mathcal B(Y, Z){\bf D}_t Z +\mathbb E 
\iint\limits_{\mathcal M\times \mathcal N}\mathcal A(Y)\mathcal B(Y,Z)+\frac{1}{2}\mathbb E\iint\limits_{\mathcal M 
\times \mathcal N} |\mathcal B(Y,Z)|^2.
\end{align}

\vspace{2mm}

{\em Step 2.  Estimates for the terms appearing  in $\mathbb E\iint\limits_{\mathcal M\times \mathcal N}\mathcal B(Y, Z){\bf D}_t Z$}

We set
\begin{align}\label{3.12}
    &\mathbb E\iint\limits_{\mathcal M\times \mathcal N}\mathcal B(Y, Z){\bf D}_t Z=K_{1}+K_{2}+K_{3},
\end{align}
where $K_i$ is the expectation of the inner product of $\mathcal B_i$ and ${\bf D}_t Z$.

For the first term $K_{1}$, we apply identity (\ref{2.10}) for time difference operator and discrete integration by parts formula (\ref{2.12}) to obtain
\begin{align}\label{3.13}
    K_1=\ & \mathbb E\iint\limits_{\mathcal M\times\mathcal N} r{\bf A}_t{\bf D}_t\rho
\left({\bf D}_t  (Z^2)-\Delta t({\bf D}_t Z)^2\right)= B_1+D_{1}+D_2,
\end{align}    
where
\begin{align*}
    B_1 =\ & \mathbb E\iint\limits_{\mathcal M\times\partial\mathcal N} r{\bf A}_t{\bf D}_t\rho{\bf t}^+\left(|Z|^2\right)n_t,\\
    D_1 =\ &  -\mathbb E\iint\limits_{\mathcal M\times\mathcal N^*}{\bf D}_t( r{\bf A}_t{\bf D}_t\rho)|Z|^2,\\
    D_2 =\ & -\Delta t\mathbb E\iint\limits_{\mathcal M\times\mathcal N} r{\bf A}_t{\bf D}_t \rho|{\bf D}_t Z|^2.
\end{align*}
By discrete integration by parts formula (\ref{2.13}), discrete identity (\ref{2.9}) and ${\bf t}^{-}(\mathcal N)=\mathcal N^*$, we obtain for the second term $K_2$ that
\begin{align}\label{3.14}
    K_2 =\ &-2\mathbb E\iint\limits_{\mathcal M\times \mathcal N} {\bf t}^+\left({\bf t}^-(r{\bf A}_x{\bf D}_x\rho {\bf A}_x{\bf D}_x Y)\right){\bf D}_t Z\nonumber\\
        =\ & -2\mathbb E\iint\limits_{\mathcal M\times\partial\mathcal N}r{\bf A}_x{\bf D}_x\rho {\bf 
A}_x{\bf D}_x Y{\bf t}^+(Z)n_t+2\mathbb E\iint\limits_{\mathcal M\times \mathcal N^*}{\bf t}^-(r{\bf A}_x{\bf D}_x\rho) {\bf D}_t({\bf A}_x{\bf D}_x Y)Z+\nonumber\\
        \ &  2\mathbb E\iint\limits_{\mathcal M\times \mathcal N^*}{\bf D}_t(r{\bf A}_x{\bf D}_x\rho) {\bf t}^+({\bf 
 A}_x{\bf D}_x Y)Z.
\end{align}
Furthermore, applying discrete integration by part for the average operator ${\bf A}_x$ to the second term on the right-hand side of (\ref{3.14}) and noticing that $z=0$ on $\mathcal \partial \mathcal M\times \mathcal N^*$ due to $y=0$ on $\partial\mathcal M\times \overline {\mathcal N}$,  we obtain
\begin{align}\label{3.15}
    K_2 =\ & B_2+D_3+D_4+D_5,
\end{align}
where 
\begin{align*}
    B_2 =\ & -2\mathbb E\iint\limits_{\mathcal M\times\partial\mathcal N}r{\bf A}_x{\bf D}_x\rho {\bf 
A}_x{\bf D}_x Y{\bf t}^+(Z)n_t,\\
    D_3 =\ & 2\mathbb E\iint\limits_{\mathcal M^*\times \mathcal N^*}{\bf t}^-\left({\bf A}_x(r{\bf A}_x{\bf D}_x\rho)\right) {\bf D}_x Z{\bf A}_x Z,\\
    D_4 =\ & \frac{1}{2}(\Delta x)^2\mathbb E\iint\limits_{\mathcal M^*\times \mathcal N^*}{\bf t}^-\left({\bf D}_x(r{\bf 
A}_x{\bf D}_x\rho)\right) |{\bf D}_x Z|^2,\\
    \displaybreak[0] D_5 = \ & 2\mathbb E\iint\limits_{\mathcal M\times \mathcal N^*}{\bf D}_t(r{\bf A}_x{\bf D}_x\rho) {\bf t}^+({\bf 
 A}_x{\bf D}_x Y)Z.
\end{align*}
As for the last term $K_3$,  by using discrete integration by parts formula (\ref{2.12}) with respect to ${\bf D}_t$ we have
\begin{align}\label{3.16}
    K_3 =\ & -\mathbb E\iint\limits_{\mathcal M\times \mathcal N^*} {\bf D}_t\left(r\varphi^{-1}\left({\bf D}_t\varphi {\bf D}_t\rho-{\bf D}_x\varphi {\bf D}_x\rho \right)Y\right) Z+\nonumber\\
    \ &\mathbb E\iint\limits_{\mathcal M\times \partial\mathcal N} r\varphi^{-1}\left({\bf D}_t\varphi {\bf D}_t\rho-{\bf D}_x\varphi {\bf D}_x\rho \right)Y{\bf t}^+(Z)n_t    = B_3+D_6+D_7,
\end{align}
where 
\begin{align*}
    B_3 =\ & \mathbb E\iint\limits_{\mathcal M\times \partial\mathcal N} r\varphi^{-1}\left({\bf D}_t\varphi {\bf D}_t\rho-{\bf D}_x\varphi {\bf D}_x\rho \right)Y{\bf t}^+(Z)n_t,\\
    D_6 =\ & -\mathbb E\iint\limits_{\mathcal M\times \mathcal N^*}{\bf D}_t\left(r\varphi^{-1}\left({\bf D}_t\varphi {\bf D}_t\rho-{\bf D}_x\varphi {\bf D}_x\rho \right)\right){\bf t}^+(Y) Z,\\
D_7 =\ & -\mathbb E\iint\limits_{\mathcal M\times \mathcal N^*}{\bf t}^-\left(r\varphi^{-1}\left({\bf D}_t\varphi {\bf D}_t\rho-{\bf 
D}_x\varphi {\bf D}_x\rho \right)\right)|Z|^2.
\end{align*}

Therefore,  substituting (\ref{3.13}), (\ref{3.15}) and (\ref{3.16}) into (\ref{3.12}) and applying Lemmas B.1 and B.2,  we obtain
\begin{align}\label{1-3.18}
\ &\mathbb E\iint\limits_{\mathcal M\times \mathcal N}\mathcal B(Y, Z){\bf D}_t Z=\sum_{j=1}^3 B_j+\sum_{j=1}^7 D_j\nonumber\\
\geq \ & \frac{1}{2}\mathbb E\iint\limits_{\mathcal M\times\mathcal N}s\lambda\varphi\partial_t\phi r^2|g|^2+{\mathbb E\iint\limits_{\mathcal M\times \mathcal 
N^*}\left(s\lambda^2\varphi|\partial_t\phi|^2+s\lambda\varphi(\partial_x^2\phi+\partial_t^2\phi)\right)|Z|^2}-\nonumber\\
\displaybreak[0] &\mathbb E\iint\limits_{\mathcal M^*\times \mathcal N}s\lambda^2\varphi|\partial_x\phi|^2|{\bf D}_x Y|^2+\mathbb E\int\limits_{\mathcal{M} } s\lambda\varphi\left(\partial_t\phi-|\partial_x\phi|\right) {\bf t}^+(|Z|^{2})\bigg|_{t=0}-\nonumber\\
 &\mathbb E \int\limits_{\mathcal{M} ^* } s\lambda\varphi|\partial_x\phi| |{\bf D}_x Y|^{2}\bigg|_{t=0}-DT_1-BT_1,
\end{align} 
where $DT_1$ and $BT_1$ are specified in Lemmas B.1 and B.2.

{\em Step 3.  Estimates for the terms appearing in $\mathbb E 
\iint\limits_{\mathcal M\times \mathcal N}\mathcal A(Y)\mathcal B(Y,Z)$}

We now present an estimate for $\mathcal A(Y)\mathcal B(Y,Z)$ on $\mathcal M\times\mathcal N$. To do this, we set 
	\begin{equation}\label{3.18} 
		\mathbb E \iint\limits_{\mathcal M\times \mathcal N}\mathcal A(Y)\mathcal B(Y,Z)=\sum_{i=1}^{2}\sum_{j=1}^{3} I_{ij},
	\end{equation}
where $I_{ij}$ is the expectation of the inner product of $\mathcal A_i$ and $\mathcal B_j$.

	We first give the estimate of $I_{11}+I_{12}+I_{13}$. Recalling the definitions of $\mathcal A_1$, $\mathcal B_1, \mathcal B_2, \mathcal B_3$ and using (\ref{2.6}), we obtain
		\begin{align}\label{3.19}
		&I_{11}+I_{12}+I_{13}\nonumber\\
		\displaybreak[0]=&-\mathbb{E} \iint\limits_{\mathcal M\times\mathcal N}\left(2r^2{\bf A}_t{\bf D}_t\rho {\bf t}^-(Z)-2r^2{\bf A}_x{\bf D}_x\rho {\bf A}_x{\bf D}_xY+r^2\varphi^{-1}\left({\bf D}_t\varphi{\bf D}_t\rho-{\bf D}_x\varphi {\bf D}_x\rho \right)Y\right)\nonumber\\
\ & \times\Big(\rho+\frac{1}{2}(\Delta x)^2{\bf D}_x^2\rho \Big){\bf D}_x^2Y\nonumber\\
		=&\mathbb{E} \iint\limits_{\mathcal M\times\mathcal N}\left(2q^{(11)}{\bf t}^-(Z)+2q^{(12)}{\bf A}_x{\bf D}_xY+q^{(13)}Y\right){\bf D}_x^2Y
	\end{align}
with
	\begin{align*}
		&q^{(11)}=-r{\bf A}_t{\bf D}_t\rho-\frac{1}{2}(\Delta x)^2(r{\bf A}_t{\bf D}_t\rho)(r{\bf D}_x^2\rho),\nonumber\\
		&q^{(12)}=r{\bf A}_x{\bf D}_x\rho +\frac{1}{2}(\Delta x)^2 (r{\bf A}_x{\bf D}_x\rho)(r {\bf D}_x^2\rho),\nonumber\\
		&q^{(13)}=-r\varphi^{-1}\left({\bf D}_t\varphi{\bf D}_t\rho-{\bf D}_x\varphi {\bf D}_x\rho\right)\Big(1+\frac{1}{2}(\Delta x)^2 r{\bf D}_x^2\rho \Big).
	\end{align*}
Using integration by parts (\ref{2.8}) with respect to ${\bf D}_x$ and noting that  ${\bf t}^-(Z)=0$ on $\partial\mathcal M\times \mathcal N$, we find that
	\begin{align}\label{3.20}
		2\mathbb{E} \iint\limits_{\mathcal M\times\mathcal N}q^{(11)}{\bf D}_x^2Y {\bf t}^-(Z)
		=-2\mathbb{E} \iint\limits_{\mathcal M^*\times\mathcal N} {\bf D}_x\left(q^{(11)}{\bf t}^-(Z)\right) {\bf D}_xY=D_8+D_9,
	\end{align}
where
\begin{align*}
&D_8=-2\mathbb{E} \iint\limits_{\mathcal M^*\times\mathcal N} {\bf D}_xq^{(11)}{\bf D}_x Y{\bf t}^-({\bf A}_xZ),\\
&D_9=-2\mathbb{E} \iint\limits_{\mathcal M^*\times\mathcal N} {\bf A}_xq^{(11)}{\bf D}_x Y{\bf t}^-({\bf D}_xZ).
\end{align*}
For the second term related to $q^{(12)}$ in (\ref{3.19}), we apply  $2{\bf A}_x{\bf D}_xY{\bf D}_x^2Y={\bf D}_x\left(|{\bf D}_xY|^2\right)$ and discrete integration by parts with respect to ${\bf D}_x$ to obtain
	\begin{align}\label{3.21}
		&2\mathbb{E} \iint\limits_{\mathcal M\times \mathcal N}q^{(12)}{\bf A}_x{\bf D}_x Y{\bf D}_x^2Y=\mathbb{E} \iint\limits_{\mathcal M\times \mathcal N}q^{(12)}{\bf D}_x\left(|{\bf D}_x Y|^2\right)=B_4+D_{10},
	\end{align}
where 
\begin{align*}
		   B_4= \ & \mathbb E\iint\limits_{\partial\mathcal M\times \mathcal N} q^{(12)}{\bf tr}\left(|{\bf D}_xY|^2\right)n_x,\\
		D_{10}= \ & -\mathbb{E} \iint\limits_{ \mathcal M^*\times\mathcal N}{\bf D}_xq^{(12)}|{\bf D}_xY|^2.
\end{align*}
Applying discrete integration by parts with respect to ${\bf D}_x$ twice and noting that $Y=0$ on $\partial \mathcal M\times \mathcal N$, we obtain
\begin{align}\label{3.22}
		\mathbb{E} \iint\limits_{\mathcal M\times\mathcal N}q^{(13)}Y{\bf D}_x^2Y=\ &-\frac{1}{2}\mathbb{E} \iint\limits_{\mathcal M^*\times \mathcal N}{\bf D}_x q^{(13)}{\bf D}_x\left(|Y|^2\right)-\mathbb{E} \iint\limits_{\mathcal M^*\times \mathcal N}{\bf A}_xq^{(13)}|{\bf D}_xY|^2\nonumber\\
=\ & D_{11}+D_{12},
\end{align}
	where
\begin{align*}
&D_{11}= \frac{1}{2}\mathbb{E} \iint\limits_{\mathcal M\times \mathcal N}{\bf D}_x^2 q^{(13)}|Y|^2,\\
&D_{12}= -\mathbb{E} \iint\limits_{\mathcal M^*\times \mathcal N}{\bf A}_xq^{(13)}|{\bf D}_xY|^2.
\end{align*}
Therefore, substituting (\ref{3.20})-(\ref{3.22}) into (\ref{3.19}) and applying Lemmas A.3, we obtain
\begin{align}\label{3.23}
\ & I_{11}+I_{12}+I_{13}=B_4+\sum_{j=8}^{12} D_j\nonumber\\
\geq \ &  -{\mathbb E\iint\limits_{\mathcal M\times\mathcal N^*}s\lambda^2\varphi|\partial_t\phi|^2|Z|^2}+\mathbb E\iint\limits_{\mathcal M^*\times\mathcal N}\left(s\lambda^2\varphi|\partial_x\phi|^2+s\lambda\varphi(\partial_x^2\phi+\partial_t^2\phi)\right)|{\bf D}_xY|^2+\nonumber\\
\ &\mathbb E\int\limits_{\mathcal M^*}s\lambda\varphi\partial_t\phi|{\bf D}_xY|^2\bigg|_{t=0}+ \mathbb E\int\limits_{\mathcal N}s\lambda\varphi\partial_x\phi{\bf tr}(|{\bf D}_xY|^2)\bigg |_{x=0}-\nonumber\\
\ &\mathbb E\int\limits_{\mathcal N}s\lambda\varphi\partial_x\phi{\bf tr}(|{\bf D}_xY|^2)\bigg |_{x=1}-DT_2-{BT_2}-{ BT_3},
\end{align}
where $DT_2$, $BT_2$ and $BT_3$ are specified in Lemma B.3.

Next we transfer to deal with $I_{21}+I_{22}+I_{23}$. We use the definitions of $\mathcal A_2$, $\mathcal B_1, \mathcal B_2, \mathcal B_3$ to rewrite $I_{21}+I_{22}+I_{23}$ as 
	\begin{align}\label{3.24}
	\	&I_{21}+I_{22}+I_{23}\nonumber\\
		=\ &\mathbb{E} \iint\limits_{\mathcal M\times\mathcal N}\left(2r{\bf A}_t{\bf D}_t\rho {\bf t}^-(Z)-2r {\bf A}_x{\bf D}_x\rho {\bf A}_x{\bf D}_x Y+r\varphi^{-1}\left({\bf D}_t\varphi {\bf D}_t\rho-{\bf D}_x\varphi {\bf D}_x\rho \right)Y\right)\nonumber\\
\ & \times\left(r{\bf D}_t^2 \rho-r{\bf D}_x^2\rho \right)Y\nonumber\\
		=\ & \mathbb{E} \iint\limits_{\mathcal M\times\mathcal N}\left(2q^{(21)}{\bf t}^-(Z)+2q^{(22)}{\bf A}_x{\bf D}_xY+q^{(23)}Y\right)Y
	\end{align}
	with
	\begin{align*}
		&q^{(21)}=r{\bf A}_t{\bf D}_t\rho\left(r{\bf D}_t^2 \rho-r{\bf D}_x^2\rho \right),\nonumber\\
		&q^{(22)} =-r {\bf A}_x{\bf D}_x\rho \left(r{\bf D}_t^2 \rho-r{\bf D}_x^2\rho \right),\nonumber\\
		&q^{(23)}=r\varphi^{-1}\left({\bf D}_t\varphi {\bf D}_t\rho-{\bf D}_x\varphi {\bf D}_x\rho \right)\left(r{\bf D}_t^2 \rho-r{\bf D}_x^2\rho \right).
	\end{align*}
Using $2{\bf t}^+ Y{\bf D}_tY={\bf D}_t\left(|Y|^2\right)+\Delta t|{\bf D}_tY|^2$ and integration by parts with respect to ${\bf D}_t$ given by (\ref{2.12}), we obtain
\begin{align}\label{3.25}
\ &\mathbb{E} \iint\limits_{\mathcal M\times\mathcal N}2q^{(21)}Y{\bf t}^-(Z)=\mathbb{E} \iint\limits_{\mathcal M\times\mathcal N^*}2{\bf t}^+\big(q^{(21)}\big){\bf t}^+(Y) {\bf D}_t Y\nonumber\\
=\ & \mathbb{E} \iint\limits_{\mathcal M\times\mathcal N^*}{\bf t}^+\big(q^{(21)}\big){\bf D}_t(|Y|^2)+\Delta t\mathbb E\iint\limits_{\mathcal M\times \mathcal N^*} {\bf t}^+\big(q^{(21)}\big)|{\bf D}_t Y|^2=B_5+D_{13}+D_{14},
\end{align}
where
\begin{align*}
B_5=\ &{\mathbb E\iint\limits_{\mathcal M\times\partial\mathcal N}{\bf t}^+\big({\bf t}^+\big(q^{(21)}\big)\big)|Y|^2 n_t},\\
D_{13}=\ & -\mathbb E\iint\limits_{\mathcal M\times\mathcal N}{\bf t}^+\big({\bf D}_tq^{(21)}\big)|Y|^2,\\
D_{14}=\ & {\Delta t\mathbb E\iint\limits_{\mathcal M\times \mathcal N^*} {\bf t}^+\big(q^{(21)}\big)|Z|^2.}
\end{align*}
From integrations by parts with respect to ${\bf A}_x$ and ${\bf D}_x$, together with $Y=0$ on $\partial \mathcal M\times\mathcal N$, it follows that
\begin{align}\label{3.26}
\ &\mathbb{E} \iint\limits_{\mathcal M\times\mathcal N}2q^{(22)}{\bf A}_x{\bf D}_x Y Y\nonumber\\
=\ &\mathbb{E} \iint\limits_{\mathcal M^*\times\mathcal N}2{\bf A}_xq^{(22)}{\bf A}_x Y{\bf D}_xY+\frac{1}{2}(\Delta x)^2\mathbb{E} \iint\limits_{\mathcal M^*\times\mathcal N}{\bf D}_xq^{(22)}|{\bf D}_x Y|^2=D_{15}+D_{16},  \nonumber\\
\end{align}
where
\begin{align*}
D_{15}=\ &-\mathbb E\iint\limits_{\mathcal M\times\mathcal N}{\bf A}_x{\bf D}_x q^{(22)}|Y|^2,\\
D_{16}=\ & \frac{1}{2}(\Delta x)^2\mathbb{E} \iint\limits_{\mathcal M^*\times\mathcal N}{\bf D}_xq^{(22)}|{\bf D}_x Y|^2.
\end{align*}
Additionally, we set
\begin{align}\label{3.27}
D_{17}=\mathbb{E} \iint\limits_{\mathcal M\times\mathcal N}q^{(23)}|Y|^2.
\end{align}
Therefore, combining the estimates (\ref{3.25})-(\ref{3.27}) and (\ref{3.24}) and applying Lemma B.4, we obtain
\begin{align}\label{3.28}
&I_{21}+I_{22}+I_{23}=B_5+\sum_{j=13}^{17} D_j\nonumber\\
\geq \ & c_0\mathbb E\iint\limits_{\mathcal M\times\mathcal N}s^3\lambda^3\varphi^3|Y|^2+\mathbb E\int\limits_{\mathcal M}\left(s^3\lambda^3\varphi^3\partial_t\phi(|\partial_t\phi|^2-|\partial_x\phi|^2)\right)|Y|^2\bigg |_{t=0}-DT_3-{ BT_4},
\end{align}
where $c_0$, $DT_3$ and $BT_4$ are specified in Lemma B.4.

Consequently, substituting (\ref{3.23}) and (\ref{3.28}) into (\ref{3.18}) yields
\begin{align}\label{3.30} 
	\ &\mathbb E \iint\limits_{\mathcal M\times \mathcal N}\mathcal A(Y)\mathcal B(Y,Z)=\sum_{j=4}^{5}B_j+\sum_{j=8}^{17}D_j\nonumber\\
\displaybreak[0]\geq \ &  c_0\mathbb E\iint\limits_{\mathcal M\times\mathcal N}s^3\lambda^3\varphi^3|Y|^2-{\mathbb E\iint\limits_{\mathcal M\times\mathcal N^*}s\lambda^2\varphi|\partial_t\phi|^2|Z|^2}+\nonumber\\
\ &\mathbb E\iint\limits_{\mathcal M^*\times\mathcal N}\left(s\lambda^2\varphi|\partial_x\phi|^2+s\lambda\varphi(\partial_x^2\phi+\partial_t^2\phi)\right)|{\bf D}_xY|^2+\nonumber\\
\ &\mathbb E\int\limits_{\mathcal M}\left(s^3\lambda^3\varphi^3\partial_t\phi(|\partial_t\phi|^2-|\partial_x\phi|^2)\right)|Y|^2\bigg |_{t=0}+\mathbb E\int\limits_{\mathcal M^*}s\lambda\varphi\partial_t\phi|{\bf D}_xY|^2\bigg|_{t=0}+\nonumber\\
\ & \mathbb E\int\limits_{\mathcal N}s\lambda\varphi\partial_x\phi{\bf tr}(|{\bf D}_xY|^2)\bigg |_{x=0}-\mathbb E\int\limits_{\mathcal N}s\lambda\varphi\partial_x\phi{\bf tr}(|{\bf D}_xY|^2)\bigg |_{x=1}-\nonumber\\
\ &DT_2-DT_3-{ BT_2}-{ BT_3}-{ BT_4}.
	\end{align}

{\em Step 4.  The remainder of the proof.}

	Obviously, for $\phi$ we have 
\begin{align*}
\partial_x\phi=2(x-x^*)<0,\quad \partial_t\phi=2\beta{(T-t)}\geq 0,\quad \partial_x^2\phi+\partial_t^2\phi=2(1-\beta)>0.
\end{align*}
By (\ref{1-3.2}), we further obtain
\begin{align*}
&\left(\partial_t\phi-|\partial_x\phi|\right) \big|_{t=0}\geq 2(\beta { T}-\sup_{x\in I}|x-x^*|)>0,\\
&\partial_t\phi\left(|\partial_t\phi|^2-|\partial_x\phi|^2\right) \big|_{t=0}\geq 8\beta  {T}(\beta^2{T^2}-\sup_{x\in I}|x-x^*|^2)>0.
\end{align*}
 Then we add up (\ref{1-3.18}) and (\ref{3.30}) to obtain 
	\begin{align}\label{3.31}
	\ &C\mathbb E\iint\limits_{\mathcal M\times \mathcal N}\mathcal B(Y, Z){\bf D}_t Z+C\mathbb E \iint\limits_{\mathcal M\times \mathcal N}\mathcal A(Y)\mathcal B(Y,Z)\nonumber\\
		\ge	\	& \mathbb{E} \iint\limits_{\mathcal M\times\mathcal N}s^3\lambda^3\varphi^3|Y|^2+{\mathbb{E} \iint\limits_{\mathcal M\times\mathcal N^*}s\lambda\varphi |Z|^2}+\mathbb{E} \iint\limits_{\mathcal M^*\times\mathcal N}s\lambda\varphi |{\bf D}_xY|^2+ \nonumber\\
		\displaybreak[0]&\mathbb{E} \iint\limits_{\mathcal M\times\mathcal N}s\lambda\varphi {(T-t)} r^2|g|^2+\mathbb{E} \int\limits_{\mathcal M}s^3\lambda^3\varphi^3|Y|^2\bigg|_{t=0}+\mathbb{E} \int\limits_{\mathcal M}s\lambda \varphi{\bf t}^+\left(|Z|^2\right)\bigg |_{t=0}+\nonumber\\
& \mathbb{E} \int\limits_{\mathcal M^*}s\lambda \varphi |{\bf D}_xY|^2\bigg |_{t=0}+\mathbb E\int\limits_{\mathcal N}s\lambda\varphi{\bf tr}(|{\bf D}_xY|^2)\bigg |_{x=1}-\mathbb{E} \int\limits_{\mathcal N}s\lambda\varphi {\bf tr}(|{\bf D}_xY|^2)\bigg|_{x=0}-\nonumber\\
\ &\sum_{j=1}^3 DT_j-\sum_{j=1}^4 BT_j.
	\end{align}
 Recalling the expressions of $DT_j$ and $BT_j$ and choosing $s\geq \mathcal O_\lambda(1)$ and $\varepsilon \mathcal O_{\lambda}(1)<1$, we obtain
\begin{align}\label{3.32}
\sum_{j=1}^3 DT_j \leq\ & Cs^3\mathbb{E} \iint\limits_{\mathcal M\times\mathcal N}|Y|^2+{Cs \iint\limits_{\mathcal M\times\mathcal N^*}|Z|^2}+Cs \iint\limits_{\mathcal M^*\times\mathcal N}|{\bf D}_xY|^2+\nonumber\\
\ &C(\Delta x)^2\mathbb E\iint\limits_{\mathcal M^*
\times\mathcal N^*}s\lambda^2\varphi|{\bf D}_x Z|^2+{Cs\Delta t}\mathbb E\iint\limits_{\mathcal M\times\mathcal N}r^2|g|^2+\nonumber\\
\ &\varepsilon\mathcal O_{\lambda}(s(\Delta x)^2) \mathbb E\iint\limits_{\mathcal M\times\mathcal N}\left(r^2|f|^2+|\mathcal 
A(Y)|^2+|{\mathcal B}(Y, Z)|^2+|\mathcal R(Y,Z)|^2\right)
\end{align}
and
\begin{align}\label{3.33}
\sum_{j=1}^4 BT_j \leq\ & Cs^3\mathbb E\int\limits_{\mathcal M}|Y|^2\bigg |_{t=0}+Cs\mathbb E\int\limits_{\mathcal M}{\bf t}^+\left(|Z|^2\right)\bigg |_{t=0}+Cs\mathbb E\int\limits_{\mathcal M^*}|{\bf D}_xY|^2\bigg |_{t=0}+\nonumber\\
\ & {C(\lambda)s^2e^{C(\lambda) s}\bigg(\mathbb E\int\limits_{\mathcal M^*}|{\bf D}_x Y|^2\bigg|_{t=T}\bigg)^{\frac{1}{2}}+C(\lambda) s^3e^{C(\lambda) s}\bigg(\mathbb E\int\limits_{\mathcal M}|Y|^2\bigg|_{t=T}\bigg)^{\frac{1}{2}}+}\nonumber\\
\ &Cs\mathbb E\int\limits_{\mathcal N}{\bf tr}(|{\bf D}_xY|^2)\bigg|_{x=0}+Cs\mathbb E\int\limits_{\mathcal N}{\bf tr}(|{\bf D}_xY|^2)\bigg|_{x=1}.
\end{align}
Therefore, substituting (\ref{3.31}) into (\ref{1-3.12}) and choosing $\lambda$ sufficiently large to absorb some terms on the right-hand side of (\ref{3.32}) and (\ref{3.33}), we find that
	\begin{align}\label{3.34}
		\hspace{-1cm} \ &\mathbb{E} \iint\limits_{\mathcal M\times\mathcal N}s^3\lambda^3\varphi^3|Y|^2+{\mathbb{E} \iint\limits_{\mathcal M\times\mathcal N^*}s\lambda\varphi |Z|^2}+\mathbb{E} \iint\limits_{\mathcal M^*\times\mathcal N}s\lambda\varphi |{\bf D}_xY|^2+ \nonumber\\
		\hspace{-1cm} \ &\mathbb{E} \iint\limits_{\mathcal M\times\mathcal N}s\lambda\varphi r^2{(T-t)}|g|^2+\mathbb{E} \int\limits_{\mathcal M}s^3\lambda^3\varphi^3|Y|^2\bigg|_{t=0}+\mathbb{E} \int\limits_{\mathcal M}s\lambda \varphi{\bf t}^+\left(|Z|^2\right)\bigg |_{t=0}+\nonumber\\
\hspace{-1cm} & \mathbb{E} \int\limits_{\mathcal M^*}s\lambda \varphi |{\bf D}_xY|^2\bigg |_{t=0}+\mathbb E\int\limits_{\mathcal N}s\lambda\varphi{\bf tr}(|{\bf D}_xY|^2)\bigg |_{x=1}\nonumber\\
\leq \ & C \mathbb E\iint\limits_{\mathcal M\times\mathcal N}r^2|f|^2+ \varepsilon\mathcal O_{\lambda}(s(\Delta x)^2) \mathbb E\iint\limits_{\mathcal M\times\mathcal N}\left(|\mathcal 
A(Y)|^2+|{\mathcal B}(Y, Z)|^2\right)+C \mathbb E\iint\limits_{\mathcal M\times\mathcal N}|\mathcal R(Y,Z)|^2+\nonumber\\
\ & C(\Delta x)^2\mathbb E\iint\limits_{\mathcal M^*
\times\mathcal N^*}s\lambda^2\varphi|{\bf D}_x Z|^2+C\mathbb E\int\limits_{\mathcal N}s\lambda\varphi{\bf tr}(|{\bf D}_xY|^2)\bigg|_{x=0}+{C(\Delta t)^{\frac{1}{2}}\mathbb E\int\limits_{\mathcal M}r^2|g|^2\bigg|_{t=T}}+\nonumber\\
\ &  {C(\lambda)s^2e^{C(\lambda) s}\bigg(\mathbb E\int\limits_{\mathcal M^*}|{\bf D}_x Y|^2\bigg|_{t=T}\bigg)^{\frac{1}{2}}+C(\lambda) s^3e^{C(\lambda) s}\bigg(\mathbb E\int\limits_{\mathcal M}|Y|^2\bigg|_{t=T}\bigg)^{\frac{1}{2}}},
	\end{align}	
{ where we have used
\begin{align*}
Cs\Delta t\mathbb E\iint\limits_{\mathcal M\times\mathcal N}r^2|g|^2\leq C\mathbb E\iint\limits_{\mathcal M\times\mathcal N}s r^2(T-t)|g|^2+C(\Delta t)^{\frac{1}{2}}\mathbb E\int\limits_{\mathcal M}r^2|g|^2\bigg|_{t=T}
\end{align*}
due to $s(\Delta t)^{\frac{1}{2}}\leq \varepsilon^{\frac{1}{2}} s\mathcal O(\Delta x)\leq C$.
}

Next, by the definitions of $\mathcal A(Y)$ and $\mathcal B(Y, Z)$, we have 
\begin{align}\label{3.35}
\mathbb E\iint\limits_{\mathcal M\times\mathcal N}|\mathcal A(Y)|^2\leq\ & C\mathbb E\iint\limits_{\mathcal M\times\mathcal N}\left(1+\mathcal O_{\lambda}((s\Delta x)^2)\right)|{\bf D}_x^2Y|^2+s^4\mathcal O_{\lambda}(1)\mathbb E\iint\limits_{\mathcal M\times\mathcal N}|Y|^2\nonumber\\
\leq\ & C(\Delta x)^{-2}\mathbb E\iint\limits_{\mathcal M^*\times\mathcal N}|{\bf D}_xY|^2+s^4\mathcal O_{\lambda}(1)\mathbb E\iint\limits_{\mathcal M\times\mathcal N}|Y|^2,
\end{align}	
and
\begin{align}\label{3.36}
\displaybreak[0]\ &\mathbb E\iint\limits_{\mathcal M\times\mathcal N}|\mathcal B(Y, Z)|^2\nonumber\\
\leq\ & s^2\mathcal O_{\lambda}(1)\mathbb E\iint\limits_{\mathcal M\times\mathcal N}{\bf t}^-\left(|Z|^2\right)+s^2\mathcal O_{\lambda}(1)\mathbb E\iint\limits_{\mathcal M\times\mathcal N}{\bf A}_x\left(|{\bf D}_xY|^2\right)+ s^2\mathcal O_{\lambda}(1)\mathbb E\iint\limits_{\mathcal M\times\mathcal N}|Y|^2\nonumber\\
\leq\ & s^2\mathcal O_{\lambda}(1)\mathbb E\int\limits_{\mathcal M}{\bf t}^+\left(|Z|^2\right)\bigg |_{t=0}+s^2\mathcal O_{\lambda}(1)\mathbb E\iint\limits_{\mathcal M\times\mathcal N}{\bf t}^+\left(|Z|^2\right)+s^2\mathcal O_{\lambda}(1)\mathbb E\iint\limits_{\mathcal M^*\times\mathcal N}|{\bf D}_xY|^2-\nonumber\\
\ &\frac{1}{2}s\mathcal O_{\lambda}(s\Delta x)\mathbb E\iint\limits_{\partial\mathcal M\times\mathcal N}{\bf tr}\left(|{\bf D}_xY|^2\right)+ s^2\mathcal O_{\lambda}(1)\mathbb E\iint\limits_{\mathcal M\times\mathcal N}|Y|^2.
\end{align}
Additionally, using the definition of  $\mathcal R(Y, Z)$ { and $r{\bf t}^+({\bf D}_t\rho)=-s\lambda\varphi(T-t)+s^2\mathcal O_{\lambda}(1)\Delta t$}, we obtain 
\begin{align}\label{3.37}
\mathbb E\iint\limits_{\mathcal M\times\mathcal N}|\mathcal R(Y,Z)|^2\leq\ & s^2\mathcal O_{\lambda}(1) \mathbb E\iint\limits_{\mathcal M\times \mathcal N}|Y|^2+{\mathbb E\iint\limits_{\mathcal M\times\mathcal N}(s^2\lambda^2\varphi^2(T-t)^2+\varepsilon)|{\bf d}_t Z|^2}.
\end{align}
On the other hand, by using (\ref{3.8}) we have
\begin{align}\label{3.38}
\ &\mathbb E\iint\limits_{\mathcal M\times\mathcal N}|{\bf d}_t Z|^2\nonumber\\
\leq \ & (\Delta t)^2 \mathbb E\iint\limits_{\mathcal M\times \mathcal N}\left(r^2|f|^2+|\mathcal 
A(Y)|^2+|{\mathcal B}(Y, Z)|^2+|\mathcal R(Y,Z)|^2\right)+\Delta t \mathbb E\iint\limits_{\mathcal M\times \mathcal N}r^2|g|^2.
\end{align}
By $s\Delta x\leq \varepsilon$ and $\Delta t=\varepsilon \mathcal O ((\Delta x)^2)$, we can choose $\varepsilon$ sufficiently small such that  
\begin{align*}
\mathcal O_{\lambda}(s^2(\Delta t)^2) \leq \varepsilon  \mathcal O_{\lambda}((\Delta x)^2)\leq\varepsilon .
\end{align*}
Then from (\ref{3.37}) and (\ref{3.38}), it follows that 
\begin{align}\label{3.39}
\mathbb E\iint\limits_{\mathcal M\times\mathcal N}|\mathcal R(Y,Z)|^2\leq \ &C \mathbb E\iint\limits_{\mathcal M\times\mathcal N}r^2|f|^2+\varepsilon\mathcal O_{\lambda}(s(\Delta x)^2) \mathbb E\iint\limits_{\mathcal M\times\mathcal N}\left(|\mathcal 
A(Y)|^2+|{\mathcal B}(Y, Z)|^2\right)+\nonumber\\
\ &  {\mathbb E\iint\limits_{\mathcal M\times\mathcal N}\varepsilon \left(s\lambda\varphi(T-t)+\Delta t\right) r^2|g|^2}+s^2\mathcal O_{\lambda}(1) \mathbb E\iint\limits_{\mathcal M\times \mathcal N}|Y|^2.
\end{align}
Consequently, combining estimates (\ref{3.35}), (\ref{3.36}) and (\ref{3.39})  yields
\begin{align}\label{3.40}
&\varepsilon\mathcal O_{\lambda}(s(\Delta x)^2) \mathbb E\iint\limits_{\mathcal M\times\mathcal N}\left(|\mathcal 
A(Y)|^2+|{\mathcal B}(Y, Z)|^2\right)+C \mathbb E\iint\limits_{\mathcal M\times\mathcal N}|\mathcal R(Y,Z)|^2\nonumber\\
\leq \ & C \mathbb E\iint\limits_{\mathcal M\times\mathcal N}r^2|f|^2+{C\mathbb E\iint\limits_{\mathcal M\times\mathcal N}\varepsilon \left(s\lambda\varphi(T-t)+\Delta t\right) r^2|g|^2}+\mathcal O_{\lambda}(1)\mathbb E\int\limits_{\mathcal N}{\bf tr}\left(|{\bf D}_xY|^2\right)\bigg|_{x=0}+\nonumber\\
\ &\mathcal O_{\lambda}(1)\mathbb E\int\limits_{\mathcal N}{\bf tr}\left(|{\bf D}_xY|^2\right)\bigg|_{x=1}+\left(\varepsilon s^3\mathcal O_{\lambda}(1)+s^2\mathcal O_{\lambda}(1)\right) \mathbb E\iint\limits_{\mathcal M\times \mathcal N}|Y|^2+\nonumber\\
\ & \varepsilon s\mathcal O_{\lambda}(1)\bigg(\mathbb E\iint\limits_{\mathcal M\times\mathcal N}|{\bf D}_xY|^2+\mathbb E\iint\limits_{\mathcal M^*\times\mathcal N}{\bf t}^+\left(|Z|^2\right)+\mathbb E\int\limits_{\mathcal M}{\bf t}^+\left(|Z|^2\right)\bigg |_{t=0}\bigg).
\end{align}

Finally, substituting (\ref{3.40}) into (\ref{3.34}) and using $\varepsilon\mathcal O_{\lambda}(1)\leq 1/2$, all terms on the right-hand side of (\ref{3.40}), except for the  terms of $f$, ${\bf tr}({\bf D}_x Y)|_{x=0}$ and $g|_{t=T}$, are absorbed into the left-hand side of (\ref{3.34}). Then going back to the origin variable $y$, we obtain the desired estimate (\ref{3.2}). This completes the proof of Theorem 3.1. \hfill$\Box$

	\section{Stability for discrete inverse problem}
	\vspace{2mm}

This section is dedicated to establishing the stability result with respect to the mesh size for our discrete inverse problem by means of Carleman estimate (\ref{3.2}).

For some known function $\overline g\in L^2_\mathcal F(\Omega; L^2(\mathcal M))$ and some constant $M>0$, we introduce the admissible set
\begin{align*}
\mathcal W=\{\ &(g, y_0, y_1)\in L^2_\mathcal F(\Omega; L^2(\mathcal M\times\mathcal N))\times L^2_\mathcal F(\Omega; H^1(\mathcal M))\times L^2_\mathcal F(\Omega; L^2(\mathcal M)): \\
 \ & g(T)=\overline g, \ \mathbb P-{\rm a.s},\ \|g\|_{L^2_\mathcal F(\Omega; L^2(\mathcal M\times\mathcal N))}+\|y_0\|_{L^2_\mathcal F(\Omega; H^1(\mathcal M))}\times \|y_1\|_{L^2_\mathcal F(\Omega; L^2(\mathcal M))}\leq R\ \}.
\end{align*}
	
\begin{thm} Let $(g^{(j)}, y_0^{(j)}, y_1^{(j)})\in \mathcal W$ for $j=1, 2$, $a, b, c, d\in L^\infty_{\mathcal F}(\Omega; L^\infty(\overline {\mathcal M}\times \overline{\mathcal N})$ and let $T$ satisfy (\ref{1-3.2}). Then there exist positive constants $C$ and sufficiently small $\widehat{\Delta x}$ and $\varepsilon$ depending on $x^*, \beta, M, R, T$ such that
		\begin{align}\label{4.1}
			\ &\left\|{\sqrt{T-t}}(g^{(1)}-g^{(2)})\right\|_{L^2_{\mathcal F}(\Omega; L^2(\mathcal M))}+\left\|y_{0}^{(1)}-y_{0}^{(2)}\right\|_{L^2_{\mathcal F}(\Omega; H^1(\mathcal M))}+\left\|y_{1}^{(1)}-y_{1}^{(2)}\right\|_{L^2_{\mathcal F}(\Omega; L^2(\mathcal M))}\nonumber\\
			\le\ & C\left\|{\bf tr}({\bf D} _{x}y^{(1)})\big|_{x=0}-{\bf tr}({\bf D} _{x}y^{(2)})\big |_{x=0}\right\|_{L_{\mathcal F}^{2}(\Omega;L^2(\mathcal N))}+C\|y^{(1)}(T)-y^{(2)}(T)\|^{\frac{1}{2}}_{L^2_\mathcal F(\Omega; H^1(\mathcal M))}+\nonumber\\
\ &C\left\| \Delta x \left({\bf D}_{t}{\bf D}_xy^{(1)}-{\bf D}_{t}{\bf D}_x y^{(2)}\right)\right\|_{L^2_{\mathcal F}(\Omega;L^2(\mathcal M^*\times \mathcal N^*))}
		\end{align}
for all $0<\Delta x\leq\widehat{\Delta x}$ and $\Delta t=\varepsilon\mathcal ((\Delta x)^2)$, where $y^{(j)}$ is the solution to (\ref{1.3})  corresponding to $\big(g^{(j)}, y_{0}^{(j)}, y_{1}^{(j)}\big)$ for $j=1,2$, respectively. 
		
	\end{thm}

	\vspace{2mm}

\noindent{\bf Remark 4.1.}\ The last error term in (\ref{4.1}) arises from the discretization of hyperbolic equation. From a control-theoretic perspective, this phenomenon stems from the lack of uniform observability estimates for discrete hyperbolic equations in the absence of appropriate penalty terms [\ref{Zuazua-SIMA-2005}]. 

\vspace{2mm}

{
\noindent{\bf Remark 4.2.}\ In discrete settings, a prior condition $g^{(1)}(T)=g^{(2)}(T)$ is required for our stability. On one hand, in both the spatially semi-discrete and fully continuous settings, i.e. $\Delta t\rightarrow 0$, this condition can be removed,  since $\Delta t$ appears as a factor in the coefficient of $g|_{t=T}$ in Carleman estimate (\ref{3.2}).   On the other hand,  $g^{(1)}(T)=g^{(2)}(T)$ can not be deduced from $\left\|\sqrt{T-t}(g^{(1)}-g^{(2)})\right\|_{L^2_{\mathcal F}(\Omega; L^2(\mathcal M))}=0$ in discrete settings. In [\ref{Lü 2015}], the corresponding continuous inverse problem was investigated under the assumption $y^{(1)}(T)=y^{(2)}(T)$. {In particular, when $g$ is independent of $t$, the last term in (\ref{3.2}) can be absorbed by the term of $g$ on the left-hand side of (\ref{3.2}), provided that $\Delta t$ is chosen sufficiently small. 
Consequently, in this case, the condition $g(T)=\bar{g}$ in $\mathcal{W}$ can be omitted.} Obviously, the assumption $y^{(1)}|_{t=T}=y^{(2)}|_{t=T}$ in [\ref{Lü 2015}] is sufficient for proving the uniqueness result, but is somewhat unreasonable for stability analysis or numerical computation. Indeed, we usually adopt $\|y^{(1)}(T)-y^{(2)}(T)\|\leq \delta$ with a small error level $\delta >0$ under a suitable norm.  Therefore, our stability result is more suitable for numerical analysis and applications.

}

\vspace{2mm}
		
	\noindent{\bf Proof.} Letting $\tilde{y}=y^{(1)}-y^{(2)}$, $\tilde{g}=g^{(1)}-g^{(2)}$ and  $\tilde{y}_{i}=y_{i}^{(1)}-y_{i}^{(2)}$ for $i=0,1$, we obtain

\begin{align}\label{4.2}
		\left\{
		\begin{array}{ll}
			{\bf d}_t \tilde z-{\bf D}_x^2\tilde y\Delta t=(a \tilde y+b{\bf A}_x{\bf D}_{x}\tilde y+c {\bf t}^{-}(\tilde z))\Delta t\\
 \hspace{3.85cm}+(d\tilde y+\tilde g) {\bf t}^+({\bf d}_tB), &(x,t)\in \mathcal M\times \mathcal N,\\
\tilde z={\bf D}_t \tilde y, &(x,t)\in \mathcal M\times \mathcal N,\\
			\tilde y(0,t)=\tilde y(L,t)=0, 	&t\in  \mathcal N, \\
			\tilde y(x,0)=\tilde y_{0},\quad {\bf t}^+(\tilde z) (x,0)=\tilde y_{1},		    &x\in {\overline {\mathcal M}}.
		\end{array}
		\right.
	\end{align}
Then applying Theorem 3.1 to $\tilde{y}$ and noticing that $g^{(1)}(T)=g^{(2)}(T)$ yields
\begin{align}\label{4.3}
		\ &\mathbb{E} \iint\limits_{\mathcal M\times\mathcal N}s^3\lambda^3\varphi^3 r^2|\tilde y|^2+{\mathbb{E} \iint\limits_{\mathcal M\times\mathcal N^*}s\lambda\varphi  r^2|{\bf D}_t \tilde y|^2}+\mathbb{E} \iint\limits_{\mathcal M^*\times\mathcal N}s\lambda\varphi r^2 |{\bf D}_x\tilde y|^2+ \nonumber\\
		&\mathbb{E} \iint\limits_{\mathcal M\times\mathcal N}s\lambda\varphi {(T-t)}r^2|\tilde g|^2+\mathbb{E} \int\limits_{\mathcal M}s^3\lambda^3\varphi^3r^2|\tilde y_0|^2+\mathbb{E} \int\limits_{\mathcal M^*}s\lambda \varphi r^2|{\bf D}_x \tilde y_0|^2+\mathbb{E} \int\limits_{\mathcal M}s\lambda \varphi r^2|\tilde y_1|^2\nonumber\\
\leq \ & C \mathbb E\iint\limits_{\mathcal M\times\mathcal N}(1+s\lambda\varphi)r^2|\tilde y|^2+{ C \mathbb E\iint\limits_{\mathcal M\times\mathcal N^*}r^2|{\bf D}_t\tilde y|^2}+C \mathbb E\iint\limits_{\mathcal M\times\mathcal N}r^2{\bf A}_x\left(|{\bf D}_x \tilde y|^2\right)+\nonumber\\
\ & C\mathbb E\int\limits_{\mathcal N}s\lambda\varphi{\bf tr}(|{\bf D}_x\tilde y|^2)\bigg|_{x=0}+C(\Delta x)^2\mathbb E\iint\limits_{\mathcal M^*
\times\mathcal N^*}s\lambda^2\varphi|{\bf D}_t{\bf D}_x \tilde y|^2+\nonumber\\
\ & { C(\lambda)s^3e^{C(\lambda)s}\|\tilde y(T)\|_{L^2_\mathcal F(\Omega; H^1(\mathcal M))}}.
	\end{align}	

On the other hand, using ${\bf A}_x(r^2)=(1+\mathcal O_{\lambda}(s\Delta x))r^2$ and discrete integration by parts (\ref{2.7}), we obtain
\begin{align}\label{4.4}
\mathbb E\iint\limits_{\mathcal M\times\mathcal N}r^2{\bf A}_x\left(|{\bf D}_x \tilde y|^2\right)=\ & \mathbb E\iint\limits_{\mathcal M^*\times\mathcal N}{\bf A}_x(r^2)|{\bf D}_x \tilde y|^2-\frac{1}{2}\Delta x\mathbb E\iint\limits_{\partial\mathcal M\times\mathcal N}r^2{\bf tr}\left(|{\bf D}_x \tilde y|^2\right)\nonumber\\
\leq \ & C\mathbb E\iint\limits_{\mathcal M^*\times\mathcal N}r^2|{\bf D}_x \tilde y|^2.
\end{align}
		
	Finally, substituting (\ref{4.3}) into (\ref{4.4}) and choosing $\lambda$ sufficiently large to absorb the first three terms on the right-hand side of (\ref{4.3}), we can obtain the desired estimate (\ref{4.1}) and then complete the proof of Theorem 4.1.\hfill$\Box$

\section{Conclusion}
\vspace{2mm}
{
This paper addresses an inverse problem associated with the fully discrete stochastic hyperbolic equation, which is a discrete counterpart of the one in [\ref{Lü 2015}]. We derive a new Carleman estimate for this fully discrete hyperbolic equation and then a conditional stability for our discrete inverse problem. 

Our work expands upon existing results on discrete Carleman estimates by extending them to the fully discrete stochastic setting.  Although the overall structure of the proof resembles that of the continuous framework in [\ref{Lü 2015}], the discrete setting introduces many additional terms that cannot be directly absorbed on the left-hand side of the Carleman estimate.  In contrast to the continuous/discrete Carleman estimates in [\ref{Carreno-ACM-2023}, \ref{Lecaros-JDE-2023}, \ref{Lü 2015}, \ref{Zhao-SICON-2025}], we have to carefully decompose the terms independent of ${\bf t}^+({\bf d}_tB)$ arising from ${\bf d}_tz-{\bf D}_x^2 y\Delta t$ to participate in computation of $\mathbb E\iint_{\mathcal M\times \mathcal N}\mathcal B(Y, Z){\bf D}_t Z$. Since the discrete Itô formula will be involved, we need to use $\mathbb E({\bf t}^+({\bf d}_t B))=0$ to eliminate discrete It\^{o} integral. This is necessary for establishing Carleman estimate, even for continuous stochastic hyperbolic equation. Another important difference is that we establish the stability result  without assumption $y^{(1)}|_{t=T}=y^{(2)}|_{t=T}$, which is more advantageous for numerical applications than the one in [28]. The lack of information about $y_t$ at the final time leads to the power $\frac{1}{2}$ of the term of $y^{(1)}(T)-y^{(2)}(T)$, which also reveals the special characteristic of our inverse problem involving random effect. In other words we can expect a Lipschitz stability for fully discrete deterministic hyperbolic equation. Also, given the discrete nature of our problem, particular attention will be devoted to the dependence on the mesh parameters $\Delta x$ and $\Delta t$, as well as to their interdependence. Several challenges remain unresolved, including stronger stability (e.g., Lipschitz stability) and the numerical analysis and implementation of this inverse problem, among others.

}

\section{Appendix}

\subsection{Appendix A}
\renewcommand{\theequation}{A.\arabic{equation}}

\setcounter{equation}{0}
{In this appendix, we show the well-posedness of the fully-discrete stochastic hyperbolic equation (\ref{1.3}). 

\vspace{2mm}

 {\noindent}{\bf Lemma A.1.}\ {\em 
		Let $a, b, c, d\in L^\infty_{\mathcal F}(\Omega; L^\infty(\overline {\mathcal M}\times \overline{\mathcal N})$, and let $g\in L^2_\mathcal F(\Omega; L^2(\mathcal M\times\mathcal N))$, $y_0\in L^{2}_{\mathcal F}(\Omega;$ $ H^{1}( {\mathcal M}))$ and $y_1\in L^{2}_{\mathcal F}(\Omega; L^{2}({\mathcal M}))$. Then there exists a unique solution $y\in L_{\mathcal{F}}^{2}(\Omega; L^\infty(\overline{\mathcal N}; H^1$ $(\mathcal M)))$ to (\ref{1.3}) such that
		\begin{align}\label{1-A.1}
			& \mathbb E \int\limits_\mathcal M {\bf t}^+(|z|)^{2}\bigg|_{t=t^n}+\mathbb E \int\limits_{\mathcal M^*}|{\bf D}_x y|^{2}\bigg|_{t=t^n}+\mathbb E \int\limits_{\mathcal M}|y|^{2}\bigg|_{t=t^n}\nonumber\\
			\le\ & C\bigg(\mathbb E \int\limits_{\mathcal M^*} |{\bf D}_xy_0|^{2}+\mathbb E \int\limits_{\mathcal M} |y_0|^{2}+\mathbb E\int\limits_{\mathcal M}| y_1|^2+\mathbb E \iint\limits_{\mathcal M\times \mathcal N} |g|^{2}\bigg),\quad n=1,2,\cdots,N
		\end{align}
for all $0<\Delta x\leq 1$ and $\Delta t=\varepsilon\mathcal O(\Delta x)$ with sufficiently small $\varepsilon$, where $C$ is positive constant depending on $a, b, c, d$ and $T$.
}

\vspace{2mm}

{\noindent\bf Proof.}\ We can rewrite the fully discrete system as
\begin{align}\label{2-A.1}
\left\{
		\begin{array}{ll}
			y^{n+1}_j=\alpha_j^n y^{n}_{j+1}+\beta_j^n y^{n}_{j}+\gamma_j^n y^{n}_{j-1}+\xi^n_j y^{n-1}_j+\zeta^n_j, &  1\leq j\leq M,\ 1\leq n\leq N,\\
			y^{n}_0=y^{n}_{M+1}=0, & 1\leq n\leq N,\\
			y^{0}_j=y_{0,j}, \quad y^1_{j}=y_{0,j}+y_{1,j}\Delta t ,  & 0\leq j\leq M+1,
		\end{array}
		\right.
\end{align}
where
\begin{align*}
&\alpha_j^n=\bigg(\frac{\Delta t}{\Delta x}\bigg)^2+\frac{b_j^n}{2\Delta x}(\Delta t)^2,\quad \beta_j^n=2-2\bigg(\frac{\Delta t}{\Delta x}\bigg)^2+\bigg(a^n_j+\frac{c^n_j}{\Delta t}\bigg)(\Delta t)^2+d^n_j(B^{n+1}-B^n)\Delta t,\\
&\gamma^n_j=\bigg(\frac{\Delta t}{\Delta x}\bigg)^2-\frac{b_j^n}{2\Delta x}(\Delta t)^2,\quad \xi^n_j=-1-c^n_j\Delta t,\quad \zeta^n_j= g^n_j(B^{n+1}-B^n)\Delta t.
\end{align*}
The existence and uniqueness of the solution to (\ref{1.3}) follow readily by recurrence. Moreover, from (\ref{2-A.1}) one can see that the value of  $y^n$ depends on $a^{n-1}_j, b^{n-1}_j, c^{n-1}_j, d^{n-1}_j, g^{n-1}_j$, $y^{n-1}$, $y^{n-2}$ and $B^{n}-B^{n-1}$. Consequently, $y^n$ is independent of $B^{n+1}-B^{n}$, which implies
\begin{align}\label{2-A.3}
\mathbb E\int\limits_{\mathcal M} {\bf t}^-(z)(dy+g){\bf t}^+({\bf d} B)=0.
\end{align}

Now we prove the discrete energy estimate (\ref{1-A.1}). Obviously, we have
\begin{align}\label{2-A.4}
{\bf d}_t(|z|^2)=2{\bf t}^-(z){\bf d}_t z+|{\bf d}_tz|^2.
\end{align}
Then multiplying the equation in (\ref{1.3}) by ${\bf t}^-(z)$ and using (\ref{2-A.3}), (\ref{2-A.4}) and integration by parts with respect to ${\bf D}_x$, we find that
\begin{align}\label{2-A.5}
&\mathbb E\int\limits_{\mathcal M}{\bf d}_t(|z|^2)+2\mathbb E\int\limits_{\mathcal M^*}{\bf t}^-({\bf D}_x z){\bf D}_x y\Delta t\nonumber\\
=\ &2\mathbb E\int\limits_{\mathcal M}{\bf t}^-(z)(a y+b {\bf A}_x{\bf D}_x y)\Delta t+2\mathbb E\int\limits_{\mathcal M} c{\bf t}^-(|z|^2)\Delta t+\mathbb E\int\limits_{\mathcal M}|{\bf d}_t z|^2\nonumber\\
\leq \ & C\Delta t\mathbb E\int\limits_{\mathcal M} {\bf t}^-(|z|^2)+C\Delta t\mathbb E\int\limits_{\mathcal M} |y|^2+C\Delta t\mathbb E\int\limits_{\mathcal M^*} |{\bf D}_x y|^2+\mathbb E\int\limits_{\mathcal M}|{\bf d}_t z|^2.
\end{align}
A direct calculation gives
\begin{align}
2\mathbb E\int\limits_{\mathcal M^*}{\bf t}^-({\bf D}_x z){\bf D}_x y\Delta t=\mathbb E\int\limits_{\mathcal M^*}{\bf t}^-\big({\bf d}_t(|{\bf D}_x y|^2)\big)+\mathbb E\int\limits_{\mathcal M^*}{\bf t}^-\big(|{\bf d}_t{\bf D}_x y|^2\big).
\end{align}
By the equation in (\ref{1.3}), we have
\begin{align}
\ &\mathbb E\int\limits_{\mathcal M}|{\bf d}_t z|^2\nonumber\\
\leq \ & C(\Delta t)^2\bigg(\mathbb E\int\limits_{\mathcal M}|{\bf D}_x^2 y|^2+\mathbb E\int\limits_{\mathcal M}|y|^2+\mathbb E\int\limits_{\mathcal M}{\bf t}^-(|z|^2)+\mathbb E\int\limits_{\mathcal M^*}|{\bf D}_x y|^2\bigg)+\nonumber\\
\ &C\Delta t\bigg(\mathbb E\int\limits_{\mathcal M}|y|^2+\mathbb E\int\limits_{\mathcal M}|g|^2\bigg)\nonumber\\
\leq \ &C\big(\varepsilon+(\Delta t)^2\big)\mathbb E\int\limits_{\mathcal M^*}|{\bf D}_x y|^2+C\Delta t\bigg(\mathbb E\int\limits_{\mathcal M}|y|^2+\mathbb E\int\limits_{\mathcal M}{\bf t}^-(|z|^2)+\mathbb E\int\limits_{\mathcal M}|g|^2\bigg).
\end{align}
Moreover, since $y|_{t=t^n}=y_0+\Delta t\sum_{k=1}^n{\bf t}^-(z)|_{t=t^k}$, we obtain
\begin{align}\label{2-A.8}
\mathbb E\int\limits_{\mathcal M}|y|^2\bigg|_{t=t^{n}}\leq C\mathbb E\int\limits_{\mathcal M}|y_0|^2+C\Delta t\sum_{k=1}^n\mathbb E\int\limits_{\mathcal M}{\bf t}^-(|z|^2)\bigg|_{t=t^k}.
\end{align}
Then, combining  (\ref{2-A.5})-(\ref{2-A.8}), we obtain 
\begin{align}
&\mathbb E\int\limits_{\mathcal M}{\bf t}^-(|z|^2)\bigg|_{t=t^{n+1}}+\mathbb E\int\limits_{\mathcal M^*}|{\bf D}_x y|^2\bigg|_{t=t^{n}}\nonumber\\
\leq\ &(1+C\Delta t)\bigg(\mathbb E\int\limits_{\mathcal M}{\bf t}^-(|z|^2)\bigg|_{t=t^{n}}+\mathbb E\int\limits_{\mathcal M^*}|{\bf D}_x y|^2\bigg|_{t=t^{n-1}}\bigg)+\nonumber\\
\ &C\Delta t\bigg(\mathbb E\int\limits_{\mathcal M}|y_0|^2+\mathbb E\int\limits_{\mathcal M}|g|^2\bigg|_{t=t^n}\bigg)+C(\Delta t)^2\sum_{k=1}^n\mathbb E\int\limits_{\mathcal M}{\bf t}^-(|z|^2)\bigg|_{t=t^k}
\end{align}
for sufficiently small $\varepsilon$. By recurrence, we derive that
\begin{align}
&\mathbb E\int\limits_{\mathcal M}{\bf t}^-(|z|^2)\bigg|_{t=t^{n+1}}+\mathbb E\int\limits_{\mathcal M^*}|{\bf D}_x y|^2\bigg|_{t=t^{n}}\nonumber\\
\leq \ &e^{CT}\bigg(\mathbb E\int\limits_{\mathcal M^*}|{\bf D}_x y_0|^2+\mathbb E\int\limits_{\mathcal M}| y_0|^2+\mathbb E\int\limits_{\mathcal M}| y_1|^2\bigg)+C\mathbb E\iint\limits_{\mathcal M\times \mathcal N}|g|^2+\nonumber\\
\ &C\Delta t\sum_{k=1}^n\mathbb E\int\limits_{\mathcal M}{\bf t}^-(|z|^2)\bigg|_{t=t^k}.
\end{align}
From discrete Gronwall's inequality [\ref{Clark-DAM-1987}], it follows that
\begin{align}
&\mathbb E\int\limits_{\mathcal M}{\bf t}^-(|z|^2)\bigg|_{t=t^{n+1}}+\mathbb E\int\limits_{\mathcal M^*}|{\bf D}_x y|^2\bigg|_{t=t^{n}}\nonumber\\
\leq \ &\prod_{k=1}^{n}(1+C\Delta t) e^{CT}\bigg(\mathbb E\int\limits_{\mathcal M^*}|{\bf D}_x y_0|^2+\mathbb E\int\limits_{\mathcal M}| y_0|^2+\mathbb E\int\limits_{\mathcal M}| y_1|^2+\mathbb E\iint\limits_{\mathcal M\times \mathcal N}|g|^2\bigg),
\end{align}
which implies the desired estimate (\ref{1-A.1}). \hfill$\Box$
}

\subsection{Appendix B}
\renewcommand{\theequation}{B.\arabic{equation}}

\setcounter{equation}{0}
 In this appendix, we first present two lemmas to show the  estimates for  $\sum_{j=1}^3B_j$ on $\mathcal M\times\partial \mathcal N$ and $\sum_{j=1}^7 D_j$ on $\mathcal M\times \mathcal N$ involved in $\mathbb E\iint\limits_{\mathcal M\times \mathcal N}\mathcal B(Y, Z){\bf D}_t Z$.

 \vspace{2mm}
 
 {\noindent}{\bf Lemma B.1.}\ {\em Provided $s\Delta x\le \varepsilon$ and $\Delta t=\varepsilon\mathcal O\left((\Delta x)^2\right)$ for a sufficiently small $\varepsilon>0$, we have
	\begin{align}\label{A.1}
\sum_{j=1}^3 B_j
\ge\ & \mathbb E\int\limits_{\mathcal{M} } s\lambda\varphi\left(\partial_t\phi-|\partial_x\phi|\right) {\bf t}^+(|Z|^{2})\bigg|_{t=0}-\mathbb E\int\limits_{\mathcal{M} ^* } s\lambda\varphi|\partial_x\phi| |{\bf D}_x Y|^{2}\bigg|_{t=0}-BT_1,
\end{align}
where
\begin{align*}
BT_1 =\ & s^2\mathcal O_{\lambda}(1)\mathbb E\int\limits_\mathcal M |Y|^2\bigg|_{t=0}+(\varepsilon s+1)\mathcal O_{\lambda}(1)\mathbb E\int\limits_\mathcal M {\bf t}^+(|Z|^2)\bigg|_{t=0}+\varepsilon s\mathcal O_{\lambda}(1)\mathbb E\int\limits_{\mathcal M^*} |{\bf D}_xY|^2\bigg|_{t=0}+\\
\ & {s^2\mathcal O_\lambda(1)e^{s\mathcal O_\lambda(1)}\bigg(\mathbb E\int\limits_{\mathcal M^*}|{\bf D}_x Y|^2\bigg|_{t=T}\bigg)^{\frac{1}{2}}+s^2\mathcal 
O_{\lambda}(1)e^{s\mathcal O_\lambda(1)}\bigg(\mathbb E\int\limits_{\mathcal M}|Y|^2\bigg|_{t=T}\bigg)^{\frac{1}{2}}.}
\end{align*}
}

\vspace{2mm}

\noindent{\bf Proof.}\ By using 
\begin{align*}
r{\bf A}_t{\bf D}_t\rho=r\partial_t\rho+s\mathcal O_{\lambda}((s\Delta t)^2)=-s\lambda \varphi\partial_t\phi+ s\mathcal O_{\lambda}((s\Delta t)^2)
\end{align*}
due to Lemma 2.5, we obtain the following estimate for $B_1$:
\begin{align}\label{A.2}
    B_1 =\ & \mathbb E\int\limits_{\mathcal M}r{\bf A}_t {\bf D}_t \rho{\bf t}^+\left(|Z|^2\right)\bigg|_{t=T}-\mathbb 
E\int\limits_{\mathcal M}r{\bf A}_t {\bf D}_t \rho {\bf t}^+\left(|Z|^2\right)\bigg|_{t=0}\nonumber\\
     \geq\ &  \mathbb E\int\limits_{\mathcal M}r{\bf A}_t {\bf D}_t \rho{\bf t}^+\left(|Z|^2\right)\bigg|_{t=T}+\mathbb 
E\int\limits_{\mathcal M}\left(s\lambda\varphi\partial_t \phi-s\mathcal O_{\lambda}((s\Delta t)^2)\right){\bf t}^+\left(|Z|^2\right)\bigg|_{t=0}.
\end{align}
{ By Lemma A.5 in [\ref{Carreno-ACM-2023}], we obtain
\begin{align}\label{1-A.3}
&r{\bf A}_t{\bf D}_t\rho=r{\bf A}_t\left(\partial_t\rho+R_{{\bf D}_t}(\rho)\right)=-s\lambda\varphi\partial_t\phi+r\tilde R(\rho),\quad(x,t)\in\mathcal M\times\mathcal N,
\end{align}
with
\begin{align}\label{1-A.4}
\tilde R(\rho)=R_{{\bf A}_t}(\partial_t\rho)+{\bf A}_t\left(R_{{\bf D}_t}(\rho)\right).
\end{align}
 A direct calculation gives
\begin{align*}
\partial_t^3\rho=\rho\left(-s\lambda^3\varphi(s^2\varphi^2-3s\varphi+1)(\partial_t\phi)^3+3s\lambda^2\varphi(s\varphi-1)
\partial_t\phi\partial_t^2\phi\right),
\end{align*} 
 which implies that
\begin{align*}
R_{{\bf A}_t}(\partial_t\rho)(x,t)=\frac{1}{4}(\Delta t)^2\int_{-1}^1 (1 - |\tau|)\partial_t^{3} \rho_{0,\tau}(x,t) {\rm d}\tau,\quad(x,t)\in\mathcal M\times\mathcal N,
\end{align*}
with $\rho_{\sigma,\tau}(x,t):=\rho(x+\sigma\frac{\Delta x}{2},t+\tau\frac{\Delta t}{2})$. Obviously, $\partial_t^3\rho(x,T)=0$ and $\partial_t^3\rho(x,T+\tau)=-\partial_t^3\rho(x,T-\tau)$ for all $\tau\in \mathbb R$, i.e. $\partial_t^3\rho_{0,\tau}(x,T)=-\partial_t^3\rho_{0,-\tau}(x,T)$. Then, we further obtain
\begin{align*}
& \int_{-1}^1 (1 - |\tau|) \partial_t^{3} \rho_{0,\tau}(x,T) {\rm d}\tau\nonumber\\
=\ &\int_{0}^1 (1 - |\tau|) \partial_t^{3} \rho_{0,\tau}(x,T) {\rm d}\tau+\int_{0}^1 (1 - |\tau|)\partial_t^{3} \rho_{0,-\tau}(x,T) {\rm d}\tau=0
\end{align*} 
and then \begin{align}\label{1-A.5}R_{{\bf A}_t}(\partial_t\rho)(x,T)=0.\end{align} 
Similarly, we have 
\begin{align*}
&{\bf A}_t\left(R_{{\bf D}_t}(\rho)\right) (x,t)= \frac{1}{16} (\Delta t)^2\int_{-1}^1 (1 - |\tau|)^2{\bf A}_t \partial_t^{3} \rho_{0,\tau}(x,t){\rm d}\tau,\quad(x,t)\in\mathcal M\times\mathcal N,
\end{align*}
and then
\begin{align}\label{1-A.6}
&{\bf A}_t\left(R_{{\bf D}_t}(\rho)\right)(x,T) \nonumber\\
=\ &   \frac{1}{8} (\Delta t)^2\int_{-1}^1 (1 - |\tau|)^2\left( \partial_t^{3} \rho_{0,\tau}\left(x,T+\frac{\Delta t}{2} \right)+\partial_t^{3} \rho_{0,\tau}\left(x,T-\frac{\Delta t}{2} \right)\right){\rm d}\tau\nonumber\\
=\ & \frac{1}{8} (\Delta t)^2\int_{0}^1 (1 - |\tau|)^2\left( \partial_t^{3} \rho_{0,\tau+1}(x,T)+\partial_t^{3} \rho_{0,-\tau+1}(x,T)\right){\rm d}\tau+\nonumber\\
\ & \frac{1}{8}(\Delta t)^2 \int_{0}^1 (1 - |\tau|)^2\left( \partial_t^{3} \rho_{0,\tau-1}(x,T)+\partial_t^{3}\rho_{0,-\tau-1} (x,T)\right){\rm d}\tau= 0.
\end{align}
Combining (\ref{1-A.3})-(\ref{1-A.6}) and $\partial_t\phi(x,T)=0$, we obtain $r{\bf A}_t{\bf D}_t\rho(x,T)=0$. Therefore, we obtain 
\begin{align}\label{1-A.7}
\mathbb E\int\limits_{\mathcal M}r{\bf A}_t {\bf D}_t \rho{\bf t}^+\left(|Z|^2\right)\bigg|_{t=T}=0.
\end{align}
Substituting (\ref{1-A.7}) into (\ref{A.2}) yields
\begin{align}\label{1-A.2}
    B_1 \geq\ &  \mathbb 
E\int\limits_{\mathcal M}\left(s\lambda\varphi\partial_t \phi-s\mathcal O_{\lambda}((s\Delta t)^2)\right){\bf t}^+\left(|Z|^2\right)\bigg|_{t=0}.
\end{align}
}

Now we transfer to estimate $B_2$. Obviously, 
\begin{align*}r{\bf A}_x{\bf D}_x\rho = -s\lambda\varphi\partial_x\phi+s\mathcal O_{\lambda}((s\Delta x)^2)=s\mathcal O_{\lambda}(1), \quad(x,t)\in\mathcal M\times\mathcal N.
\end{align*} 
{Then we have
\begin{align}\label{1-A.9}
    B_2 \geq \ & -s\mathcal O_{\lambda}(1)\mathbb E\int\limits_{\mathcal M}|{\bf A}_x{\bf D}_x Y||{\bf t}^+(Z)|\bigg|_{t=T}-\mathbb 
E\int\limits_{\mathcal M^*}\left(s\lambda\varphi|\partial_x\phi|+s\mathcal O_{\lambda}\left(s\Delta 
 x\right)^2\right)|{\bf D}_xY|^2\bigg|_{t=0}-\nonumber\\
 \ &\mathbb 
E\int\limits_{\mathcal M} \left(s\lambda\varphi|\partial_x\phi|+s\mathcal O_{\lambda}\left((s\Delta x)^2\right)\right){\bf t}^+\left(|Z|^2\right)\bigg|_{t=0}.
\end{align}
 Using H\"{o}lder inequality and $|{\bf A}_x{\bf D}_x Y|^2\leq {\bf A}_x(|{\bf D}_x Y|^2)$ yields that
 \begin{align}\label{1-A.10}
   \mathbb E\int\limits_{\mathcal M}|{\bf A}_x{\bf D}_x Y||{\bf t}^+(Z)|\bigg|_{t=T} \leq \ & 
     \bigg(\mathbb E\int\limits_{\mathcal M}{\bf A}_x({\bf D}_x Y)^2\bigg|_{t=T}\bigg)^{\frac{1}{2}} \bigg(\mathbb E\int\limits_{\mathcal M}{\bf t}^+(|Z|^2)\bigg|_{t=T}\bigg)^{\frac{1}{2}}\nonumber\\
     \leq \ &  s\mathcal O_\lambda(1)e^{s\mathcal O_\lambda(1)}\bigg(\mathbb E\int\limits_{\mathcal M^*}|{\bf D}_x Y|^2\bigg|_{t=T}\bigg)^{\frac{1}{2}}.
\end{align}
In the last inequality in (\ref{1-A.10}), we have used the discrete energy estimate in Lemma A.1 to obtain 
\begin{align}\label{1-A.11}
      \mathbb E\int\limits_{\mathcal M}{\bf t}^+(|Z|^2)\bigg|_{t=T} \leq \ & s^2\mathcal O_\lambda(1)e^{s\mathcal O_\lambda(1)}\left(\mathbb E\int\limits_{\mathcal M}{\bf t}^+(|z|^2)\bigg|_{t=T}+\mathbb E\int\limits_{\mathcal M}|y|^2\bigg|_{t=T}\right)\nonumber\\
    \leq \ & s^2\mathcal O_\lambda(1)e^{s\mathcal O_\lambda(1)}\Big(\|y_0\|^2_{L^2_\mathcal F(\Omega; H^1(\mathcal M))}+\|y_1\|^2_{L^2_\mathcal F(\Omega; L^2(\mathcal M))}+\nonumber\\
    \ &\|f\|^2_{L^2_\mathcal F(\Omega; L^2(\mathcal M))}+\|g\|^2_{L^2_\mathcal F(\Omega; L^2(\mathcal M))}\Big).
\end{align} 
Therefore, substituting (\ref{1-A.10}) into (\ref{1-A.9}), we obtain
\begin{align}\label{A.3}
    B_2 \geq\ & -s^2\mathcal O_\lambda(1)e^{s\mathcal O_\lambda(1)}\bigg(\mathbb E\int\limits_{\mathcal M^*}|{\bf D}_x Y|^2\bigg|_{t=T}\bigg)^{\frac{1}{2}}-\mathbb 
E\int\limits_{\mathcal M^*}\left(s\lambda\varphi|\partial_x\phi|+s\mathcal O_{\lambda}\left(s\Delta 
 x\right)^2\right)|{\bf D}_xY|^2\bigg|_{t=0}-\nonumber\\
 \ &\mathbb 
E\int\limits_{\mathcal M} \left(s\lambda\varphi|\partial_x\phi|+s\mathcal O_{\lambda}\left((s\Delta x)^2\right)\right){\bf t}^+\left(|Z|^2\right)\bigg|_{t=0}.
\end{align}
}

For $B_3$, we notice that 
\begin{align}\label{A.4}
       &r\varphi^{-1}({\bf D}_t\varphi {\bf D}_t \rho-{\bf D}_x\varphi {\bf D}_x \rho)\nonumber\\
       =\ & (\varphi^{-1}{\bf D}_t\varphi)( r{\bf D}_t \rho)-(\varphi^{-1}{\bf D}_x\varphi)( r{\bf D}_x \rho)\nonumber\\
    =\ & (\varphi^{-1}\partial_t\varphi+(\Delta t)^2\mathcal O_{\lambda}(1))\left(r\partial_t\rho+s\mathcal O_\lambda((s\Delta t)^2)\right)-\nonumber\\
    \ &(\varphi^{-1}\partial_x\varphi+(\Delta x)^2\mathcal O_{\lambda}(1)) \left(r\partial_x\rho+s\mathcal O_\lambda((s\Delta x)^2)\right)\nonumber\\
    =\ &-s\lambda^2\varphi|\partial_t\phi|^2+s\lambda^2\varphi|\partial_x\phi|^2+s\mathcal
O_{\lambda}((s\Delta x)^2)+s\mathcal
O_{\lambda}((s\Delta t)^2)\nonumber\\
=\ & s\mathcal O_{\lambda}(1).
\end{align}
Then using H\"{o}lder inequality and (\ref{1-A.11}) again, we obtain
{
\begin{align}\label{A.5}
    B_3 \geq\ & - s\mathcal 
O_{\lambda}(1)\bigg(\mathbb E\int\limits_{\mathcal M}|Y|^2\bigg|_{t=T}\bigg)^{\frac{1}{2}} \bigg(\mathbb E\int\limits_{\mathcal M}{\bf t}^+(|Z|^2)\bigg|_{t=T}\bigg)^{\frac{1}{2}}-\nonumber\\
     \ &  s^2\mathcal O_{\lambda}(1)\mathbb E\int\limits_{\mathcal M}|Y|^2\bigg|_{t=0}-\mathcal O_{\lambda}(1)\mathbb 
E\int\limits_{\mathcal M}{\bf t}^+\left(|Z|^2\right)\bigg|_{t=0}\nonumber\\
\geq \ & -s^2\mathcal 
O_{\lambda}(1)e^{s\mathcal O_\lambda(1)}\bigg(\mathbb E\int\limits_{\mathcal M}|Y|^2\bigg|_{t=T}\bigg)^{\frac{1}{2}}-s^2\mathcal O_{\lambda}(1)\mathbb E\int\limits_{\mathcal M}|Y|^2\bigg|_{t=0}-\nonumber\\
\ &\mathcal O_{\lambda}(1)\mathbb 
E\int\limits_{\mathcal M}{\bf t}^+\left(|Z|^2\right)\bigg|_{t=0}
\end{align}
}

Finally, combining the estimates (\ref{1-A.2}), (\ref{A.3})and (\ref{A.5}),  and using $\Delta t=\varepsilon\mathcal O((\Delta x)^2)$, we obtain the desired estimate (\ref{A.1}) and then complete the proof of Lemma B.1. \hfill$\Box$

\vspace{2mm}

{\noindent}{\bf Lemma B.2.}\ {\em Provided $s\Delta x\le \varepsilon$ and $\Delta t=\varepsilon\mathcal O\left((\Delta x)^2\right)$ for a sufficiently small $\varepsilon>0$, we have
	\begin{align}\label{A.6}
\sum_{j=1}^7 D_j
\ge\ & \frac{1}{2}\mathbb E\iint\limits_{\mathcal M\times\mathcal N}s\lambda\varphi\partial_t\phi r^2|g|^2+{\mathbb E\iint\limits_{\mathcal M\times \mathcal 
N^*}\left(s\lambda^2\varphi|\partial_t\phi|^2+s\lambda\varphi(\partial_x^2\phi+\partial_t^2\phi)\right)|Z|^2}-\nonumber\\
\ &\mathbb E\iint\limits_{\mathcal M^*\times \mathcal N}s\lambda^2\varphi|\partial_x\phi|^2|{\bf D}_x Y|^2-DT_1   
\end{align}
where
\begin{align*}
DT_1=\ & s^2\mathcal O_{\lambda}(1)\mathbb E\iint\limits_{\mathcal M\times \mathcal N}|Y|^2+{(\varepsilon s+1)\mathcal O_{\lambda}(1)\mathbb E\iint\limits_{\mathcal M\times \mathcal 
N^*}|Z|^2}+\varepsilon s\mathcal O_{\lambda}(1)\mathbb E\iint\limits_{\mathcal M^*\times \mathcal N}|{\bf D}_x Y|^2+\nonumber\\
\ &\varepsilon\mathcal O_\lambda(s(\Delta x)^2)\mathbb E\iint\limits_{\mathcal M\times\mathcal N}\left(r^2|f|^2+|\mathcal 
A(Y)|^2+|{\mathcal B}(Y, Z)|^2+|\mathcal R(Y,Z)|^2\right)+\nonumber\\
\ & {\varepsilon  s\mathcal O_{\lambda}(\Delta t)}\mathbb E\iint\limits_{\mathcal M\times\mathcal N}r^2|g|^2+\mathcal O \left((\Delta x)^2\right)\mathbb E\iint\limits_{\mathcal M^*\times\mathcal N^*}s\lambda^2\varphi|{\bf D}_x Z|^2.
\end{align*}
}

{\noindent\bf Proof.}\ We easily see 
By Lemma 2.5,  we have
\begin{align}\label{A.8}
& {\bf D}_t(r{\bf A}_t{\bf D}_t \rho)=-s\lambda^2\varphi|\partial_t\phi|^2-s\lambda\varphi\partial_t^2\phi+s\mathcal O_{\lambda}((s\Delta t)^2)=s\mathcal O_{\lambda}(1).
\end{align}
Then, we obtain the following estimate for $D_1$ that
\begin{align}\label{A.9}
    D_1 \geq\ & \mathbb E\iint\limits_{\mathcal M\times \mathcal 
N^*}\left(s\lambda^2\varphi|\partial_t\phi|^2+s\lambda\varphi\partial_t^2\phi-s\mathcal O_{\lambda}(s\Delta t)\right)|Z|^2.
\end{align}

For $D_2$, we easily see that
\begin{align}\label{A.10}
D_2 = -\frac{1}{\Delta t}\mathbb E\iint\limits_{\mathcal M\times\mathcal N}r{\bf A}_t{\bf D}_t\rho|{\bf d}_t Z|^2.
\end{align}
Substituting
\begin{align}\label{A.11}
&-r{\bf A}_t{\bf D}_t\rho=s\lambda\varphi\partial_t\phi-s\mathcal O_{\lambda}((s\Delta t)^2)
\end{align}
into (\ref{A.10}), it follows that
\begin{align}\label{A.16}
D_2 \geq  \frac{1}{\Delta t}\mathbb E\iint\limits_{\mathcal M\times\mathcal N}\left(s\lambda\varphi\partial_t\phi-s\mathcal O_{\lambda}((s\Delta t)^2)\right)|{\bf d}_t Z|^2.
\end{align}
On the other hand, from (\ref{3.8}) we obtain
\begin{align}\label{A.17}
    |{\bf d}_t Z|^2 =\ & |rf\Delta t+rg{\bf t}^+\left({\bf d}_t B\right)-\mathcal A(Y)\Delta t-\mathcal B(Y,Z)\Delta t+\mathcal R(Y,Z)\Delta t |^2\nonumber\\
    \geq \ & \frac{1}{2}r^2|g|^2({\bf t}^+ ({\bf d}_t B))^2 -4r^2|f|^2(\Delta t)^2-4|\mathcal A(Y)|^2(\Delta t)^2-\nonumber\\
                        \ & 4|\mathcal B(Y, Z)|^2(\Delta t)^2-4|\mathcal R(Y,Z)|^2(\Delta t)^2.
\end{align}
Noticing that $\partial_t\phi>0$ and $\mathbb E {({\bf t}^+ ({\bf d}_t B))^2}=\Delta t$ and substituting (\ref{A.17}) into (\ref{A.16}), we find that
\begin{align}
    D_2 \geq\ & \frac{1}{2}\mathbb E\iint\limits_{\mathcal M\times\mathcal N} s\lambda\varphi\partial_t\phi r^2|g|^2-s\mathcal O_{\lambda}((s\Delta t)^2)\mathbb E\iint\limits_{\mathcal M\times\mathcal N} r^2|g|^2-\nonumber\\
            \ & \mathcal O_\lambda(s\Delta t)\mathbb E\iint\limits_{\mathcal M\times\mathcal N}\left(r^2|f|^2+|\mathcal 
A(Y)|^2+|{\mathcal B}(Y, Z)|^2+|\mathcal R(Y,Z)|^2\right).
\end{align}

We apply discrete integration  by parts formula (\ref{2.8}) with respect to the operators ${\bf D}_x$ to obtain for $D_3$ that 
\begin{align}
\hspace{-0.3cm}D_3=\ &\mathbb E\iint\limits_{\mathcal M^*\times \mathcal N^*}{\bf t}^-\left({\bf A}_x(r{\bf A}_x{\bf D}_x\rho)\right){\bf D}_x(Z^2)=-\mathbb E\iint\limits_{\mathcal M\times \mathcal N^*}{\bf t}^-({\bf A}_x{\bf D}_x(r{\bf A}_x{\bf D}_x\rho))|Z|^2.
\end{align}
Together with ${\bf t}^-\left({\bf A}_x{\bf D}_x(r{\bf A}_x{\bf D}_x\rho)\right)=\partial_x(r\partial_x\rho)+s\mathcal O_{\lambda}((s\Delta x)^2)$ due to Lemma 2.5, we further obtain
\begin{align}
D_3\geq\ &{\mathbb E\iint\limits_{\mathcal M\times \mathcal N^*}\left(s\lambda^2\varphi|\partial_x\phi|^2+s\lambda\varphi\partial_x^2\phi-s\mathcal O_{\lambda}((s\Delta x)^2)\right)|Z|^2.}
\end{align}

By using \begin{align*}
& {\bf t}^-\left({\bf D}_x(r{\bf A}_x{\bf D}_x\rho)\right)={\bf t}^-\left(\partial_x(r\partial_x\rho)\right)+s\mathcal O_{\lambda}((s\Delta x)^2)\\
=\ &-s\lambda^2\varphi|\partial_x 
\phi|^2-s\lambda\varphi\partial_x^2\phi+\mathcal O_{\lambda}(s\Delta t)+s\mathcal O_{\lambda}((s\Delta x)^2),
\end{align*} 
we obtain the following estimate for $D_4$:
\begin{align}
    D_4\geq \ & -\frac{1}{2}(\Delta x)^2\mathbb E\iint\limits_{\mathcal M^*\times \mathcal N^*}\left(s\lambda^2\varphi|\partial_x 
\phi|^2+s\lambda\varphi\partial_x^2\phi+\mathcal O_{\lambda}(s\Delta t)+s\mathcal O_{\lambda}(s\Delta x)\right)|{\bf D}_xZ|^2.
\end{align}

 Applying Lemma 2.6 and $|{\bf A}_x{\bf D}_x Y|^2\leq {\bf A}_x\left(|{\bf D}_x Y|^2\right)$, it follows that
\begin{align}
D_5 =\ & 2\mathbb E\iint\limits_{\mathcal M\times \mathcal N^*}\left(-s\lambda^2\varphi\partial_x\phi\partial_t\phi+
\mathcal O_{\lambda}(s^2\Delta x)\right){\bf t}^+({\bf A}_x{\bf D}_x Y)Z\nonumber\\
\geq \ & -\mathbb E\iint\limits_{\mathcal M\times \mathcal N^*}\left(s\lambda^2\varphi|\partial_x\phi|^2+\mathcal O_{\lambda}(s^2\Delta x)\right){\bf t}^+({\bf A}_x(|{\bf D}_x Y|^2))-\nonumber\\
& \mathbb E\iint\limits_{\mathcal M\times \mathcal N^*}\left(s\lambda^2\varphi|\partial_t\phi|^2+\mathcal O_{\lambda}(s^2\Delta x)\right) |Z|^2.
\end{align} 
 Then by using ${\bf t}^+(\mathcal N^*)=\mathcal N$ and discrete integration by parts (\ref{2.7}) and 
\begin{align*}
{\bf t}^-\left({\bf A}_x(\varphi|\partial_x\phi|^2)\right)=\ &\varphi|\partial_x\phi|^2+\mathcal O_{\lambda}((\Delta x)^2)+\mathcal O_{\lambda}(\Delta t),
\end{align*} we further obtain 
\begin{align}
\displaybreak[0]\hspace{-2.3cm}D_5 \geq \ & -\mathbb E\iint\limits_{\mathcal M^*\times \mathcal N}\left(s\lambda^2{\bf t}^-({\bf A}_x(\varphi|\partial_x\phi|^2))+s\mathcal O_{\lambda}(s\Delta x)\right)|{\bf D}_x Y|^2+\nonumber\\
\hspace{-2.3cm} \ & \frac{1}{2}\Delta x\mathbb E\iint\limits_{\partial\mathcal M\times \mathcal N}\left(s\lambda^2{\bf t}^-(\varphi|\partial_x\phi|^2)+s\mathcal O_{\lambda}(s\Delta x)\right){\bf tr}\left(|{\bf D}_x Y|^2\right)-\nonumber\\
\ & { \mathbb E\iint\limits_{\mathcal M\times \mathcal N^*}\left(s\lambda^2\varphi|\partial_t\phi|^2+s\mathcal O_{\lambda}(s\Delta x)\right)|Z|^2}\nonumber\\
\geq \ &  -\mathbb E\iint\limits_{\mathcal M^*\times \mathcal N}\left(s\lambda^2\varphi|\partial_x\phi|^2+s\mathcal O_{\lambda}(s\Delta x)\right)|{\bf D}_x Y|^2-\nonumber\\
\ & {\mathbb E\iint\limits_{\mathcal M\times \mathcal N^*}\left(s\lambda^2\varphi|\partial_t\phi|^2+s\mathcal O_{\lambda}(s\Delta x)\right)|Z|^2},
\end{align}
where we have used 
\begin{align*}
s\lambda^2{\bf t}^{-}(\varphi|\partial_x\phi|^2)+s\mathcal O_{\lambda}(s\Delta x)> 0,
\end{align*} if $s\Delta x\leq \varepsilon$ with sufficiently small $\varepsilon>0$.

Similar to (\ref{A.4}),  we have
\begin{align}\label{A.24}
 \ &{\bf D}_t(r\varphi^{-1}({\bf D}_t\varphi {\bf D}_t \rho-{\bf D}_x\varphi {\bf D}_x \rho))\nonumber\\
=\ &\partial_t\left(-s\lambda^2\varphi|\partial_t\phi|^2+s\lambda^2\varphi|\partial_x\phi|^2\right)+s\mathcal
O_{\lambda}(s\Delta x)=s\mathcal O_{\lambda}(1).
\end{align}
Then combining (\ref{A.4}) and (\ref{A.24}), we obtain the following estimates for $D_6$ and $D_7$
\begin{align}
D_6\geq \ &-s^2\mathcal O_{\lambda}(1)\mathbb E\iint\limits_{\mathcal M\times \mathcal N^*}{\bf t}^+(|Y|^2)-\mathcal O_{\lambda}(1)\mathbb E\iint\limits_{\mathcal M\times \mathcal N^*}|Z|^2\nonumber\\
\geq \ & -s^2\mathcal O_{\lambda}(1)\mathbb E\iint\limits_{\mathcal M\times \mathcal N}|Y|^2-{\mathcal O_{\lambda}(1)\mathbb E\iint\limits_{\mathcal M\times \mathcal N^*}|Z|^2}
\end{align}
and
\begin{align}
D_7
\geq \ & {\mathbb E\iint\limits_{\mathcal M\times\mathcal N^*}s\lambda^2\varphi(|\partial_t\phi|^2-|\partial_x\phi|^2)|Z|^2-s\mathcal O_{\lambda}(s\Delta x) \mathbb E\iint\limits_{\mathcal M\times\mathcal N^*}|Z|^2.}
\end{align}

Finally, combining the preceding estimates for $D_1, \cdots, D_7$ and noticing that $s\Delta x\leq \varepsilon$ and $\Delta t=\varepsilon\mathcal O((\Delta x)^2)$, we obtain the desired estimate (\ref{A.6}) and then complete the proof of Lemma B.2. \hfill$\Box$

\vspace{2mm}

We now present the estimates of boundary terms $B_4, B_5$ and distribution terms $D_8, \cdots$, $D_{17}$ appearing in $\mathbb E 
\iint\limits_{\mathcal M\times \mathcal N}\mathcal A(Y)\mathcal B(Y,Z)$.

\vspace{2mm}

{\noindent}{\bf Lemma B.3.}\ {\em Provided $s\Delta x\le \varepsilon$ and $\Delta t=\varepsilon\mathcal O\left((\Delta x)^2\right)$ for a sufficiently small $\varepsilon>0$, we have
	\begin{align}\label{A.27}
B_4\geq \ & \mathbb E\int\limits_{\mathcal N}s\lambda\varphi\partial_x\phi{\bf tr}(|{\bf D}_xY|^2)\bigg |_{x=0}-\mathbb E\int\limits_{\mathcal N}s\lambda\varphi\partial_x\phi{\bf tr}(|{\bf D}_xY|^2)\bigg |_{x=1}-{ BT_2}
\end{align}
and
\begin{align}\label{A.28}
\sum_{j=8}^{12} D_j \ge\ & -{\mathbb E\iint\limits_{\mathcal M\times\mathcal N^*}s\lambda^2\varphi|\partial_t\phi|^2|Z|^2}+\mathbb E\iint\limits_{\mathcal M^*\times\mathcal N}\left(s\lambda^2\varphi|\partial_x\phi|^2+s\lambda\varphi(\partial_x^2\phi+\partial_t^2\phi)\right)|{\bf D}_xY|^2+\nonumber\\
\ &\mathbb E\int\limits_{\mathcal M^*}s\lambda\varphi\partial_t\phi|{\bf D}_xY|^2\bigg |_{t=0}-DT_2-{ BT_3} ,
\end{align}
where
\begin{align*}
DT_2=\ & s\mathcal O_{\lambda}(1)\mathbb E\iint\limits_{\mathcal M\times \mathcal N}|Y|^2+{\varepsilon s\mathcal O_{\lambda}(1)\mathbb E\iint\limits_{\mathcal M\times \mathcal 
N^*}|Z|^2}+\varepsilon s\mathcal O_{\lambda}(1)\mathbb E\iint\limits_{\mathcal M^*\times \mathcal N}|{\bf D}_x Y|^2+\nonumber\\
\ & \mathcal O((\Delta x)^2)\mathbb E\iint\limits_{\mathcal M^*\times\mathcal N^*}s\lambda\varphi|{\bf D}_xZ|^2
\end{align*}
and
\begin{align*}
{BT_2} =\ & \varepsilon s\mathcal O_{\lambda}(1)\mathbb E\int\limits_{\mathcal N}{\bf tr}(|{\bf D}_xY|^2)\bigg |_{x=0}+\varepsilon s\mathcal O_{\lambda}(1)\mathbb E\int\limits_{\mathcal N}{\bf tr}(|{\bf D}_xY|^2)\bigg |_{x=1},\nonumber\\
 { BT_3}= \ &\varepsilon s\mathcal O_{\lambda}(1)\mathbb E\int\limits_{\mathcal M^*}|{\bf D}_xY|^2 \bigg|_{t=0}+ {s^2\mathcal O_{\lambda}(1)e^{s\mathcal O_{\lambda}(1)}\bigg(\mathbb E\int\limits_{\mathcal M^*} |{\bf D}_xY|^2\bigg|_{t=T}\bigg)^{\frac{1}{2}}}.
\end{align*}}

\vspace{2mm}

{\noindent\bf Proof.}\ By Lemma 2.5 and Lemma 2.6, we have
   \begin{align}
   {\bf D}_{x} q^{(11)}=\ &-{\bf D}_x(r{\bf A}_t{\bf D}_t\rho)-\frac{1}{2}(\Delta x)^2\left({\bf D}_x(r{\bf A}_t{\bf D}_t\rho){\bf A}_x(r{\bf D}_x^2\rho)+{\bf A}_x(r{\bf A}_t{\bf D}_t\rho){\bf D}_x(r{\bf D}_x^2\rho)\right)\nonumber\\
   =\ & -\partial_x(r\partial_t\rho)+\mathcal O_{\lambda}(s^2\Delta s)+(\Delta x)^2 s^3\mathcal O_{\lambda}(1)= s\lambda^2\varphi\partial_x\phi\partial_t\phi+s\mathcal O_{\lambda}(s\Delta x).
   \end{align}
	Then, using Young's inequality and $|{\bf t}^-({\bf A}_xZ)|^2\leq {\bf A}_x{\bf t}^-(|Z|^2)$, we obtain 
	\begin{align}\label{A.31}
		D_8		\ge \ &- \mathbb{E} \iint\limits_{\mathcal M^*\times\mathcal N}\left(s\lambda^2\varphi |\partial_x\phi|^2+s\mathcal O_{\lambda}(s\Delta x)\right)|{\bf D}_xY|^2- \nonumber\\
\ & \mathbb{E} \iint\limits_{\mathcal M^*\times\mathcal N}\left(s\lambda^2\varphi |\partial_t\phi|^2+s\mathcal O_{\lambda}(s\Delta x)\right){\bf A}_x{\bf t}^-\left(|Z|^2\right),
	\end{align}
On the other hand, applying integration by parts with respect to ${\bf A}_x$ to the second term on the right-hand side of (\ref{A.31}), and noting that ${\bf t}^-(Z)=0$ on $\partial \mathcal M\times \mathcal N$ and ${\bf t}^+({\bf A}_x\varphi|\partial_t\phi|^2)=\varphi|\partial_t\phi|^2+\mathcal O_{\lambda}(\Delta t)+\mathcal O_{\lambda}((\Delta x)^2)$,  we find that
\begin{align}\label{A.32}
&-\mathbb{E} \iint\limits_{\mathcal M^*\times\mathcal N}\left(s\lambda^2\varphi |\partial_t\phi|^2+s\mathcal O_{\lambda}(s\Delta x)\right){\bf A}_x{\bf t}^-\left(|Z|^2\right)\nonumber\\
=\ & -\mathbb{E} \iint\limits_{\mathcal M\times\mathcal N^*}\left(s\lambda^2{\bf t}^+\left({\bf A}_x\varphi |\partial_t\phi|^2\right)+s\mathcal O_{\lambda}(s\Delta x)\right)|Z|^2\nonumber\\
\geq\ &{-\mathbb{E} \iint\limits_{\mathcal M\times\mathcal N^*}\left(s\lambda^2\varphi |\partial_t\phi|^2+s\mathcal O_{\lambda}(s\Delta x)\right)|Z|^2.}
\end{align}
Therefore, combining (\ref{A.31}) and (\ref{A.32}) we obtain the following estimate for $D_8$:
\begin{align}\label{A.33}
D_8		\ge \ &- \mathbb{E} \iint\limits_{\mathcal M^*\times\mathcal N}\left(s\lambda^2\varphi |\partial_x\phi|^2+s\mathcal O_{\lambda}(s\Delta x)\right)|{\bf D}_xY|^2-\nonumber\\
\ & { \mathbb{E} \iint\limits_{\mathcal M\times\mathcal N^*}\left(s\lambda^2\varphi |\partial_t\phi|^2+s\mathcal O_{\lambda}(s\Delta x)\right)|Z|^2.}
\end{align}

For $D_9$,  applying integration by parts formulas (\ref{2.11}) and (\ref{1-2.14})  with respect to ${\bf t}^-$ and ${\bf D}_t$ and using $2{\bf t}^+({\bf D}_x Y){\bf D}_x Z={\bf D}_t\left(|{\bf D}_x Y|^2\right)+\Delta t|{\bf D}_x Z|^2$ leads to
		\begin{align}\label{}
		D_9 =\ &-\mathbb{E} \iint\limits_{\mathcal M^*\times \mathcal N^*} {\bf t}^+\big({\bf A}_x q^{(11)}\big){\bf D}_t(|{\bf D}_x Y|^2)-\Delta t\mathbb{E} \iint\limits_{\mathcal M^*\times \mathcal N^*} {\bf t}^+\big({\bf A}_x q^{(11)}\big)|{\bf D}_x Z|^2\nonumber\\
		=\ &\mathbb{E} \iint\limits_{\mathcal M^*\times \mathcal N^*} {\bf D}_t{\bf A}_x q^{(11)}{\bf t}^-(|{\bf D}_x Y|^2)-\mathbb{E} \iint\limits_{\mathcal M^*\times \partial\mathcal N} {\bf A}_x q^{(11)}|{\bf D}_x Y|^2n_t-\nonumber\\
\ & \Delta t\mathbb{E} \iint\limits_{\mathcal M^*\times \mathcal N^*} {\bf t}^+\big({\bf A}_x q^{(11)}\big)|{\bf D}_x Z|^2\nonumber\\
=\ &\mathbb{E} \iint\limits_{\mathcal M^*\times \mathcal N} {\bf t}^+\big({\bf D}_t{\bf A}_x q^{(11)}\big)|{\bf D}_x Y|^2-\Delta t\mathbb{E} \iint\limits_{\mathcal M^*\times \partial\mathcal N} {\bf t}^+\big({\bf D}_t{\bf A}_x q^{(11)}\big)|{\bf D}_x Y|^2n_t-\nonumber\\
\ & \mathbb{E} \iint\limits_{\mathcal M^*\times \partial\mathcal N} {\bf A}_x q^{(11)}|{\bf D}_x Y|^2n_t-\Delta t\mathbb{E} \iint\limits_{\mathcal M^*\times \mathcal N^*} {\bf t}^+\big({\bf A}_x q^{(11)}\big)|{\bf D}_x Z|^2.
	\end{align}
Lemma 2.5 and Lemma 2.6 gives
\begin{align*}
{\bf D}_t q^{(11)}=\ & -{\bf D}_t(r{\bf A}_t{\bf D}_t\rho)-\frac{1}{2}(\Delta x)^2\left({\bf D}_t(r{\bf A}_t{\bf D}_t\rho){\bf t}^+(r{\bf D}_x^2\rho)+{\bf t}^-(r{\bf A}_t{\bf D}_t\rho){\bf D}_t(r{\bf D}_x^2\rho)\right)\\
=\ & -\partial_t(r\partial_t\rho)+s\mathcal O_{\lambda}\left((s\Delta t)^2\right)-\frac{1}{2}(\Delta x)^2s\mathcal O_{\lambda}(1)s^2\mathcal O_{\lambda}(1)\\
=\ &s\lambda^2\varphi|\partial_t\phi|^2+s\lambda\varphi\partial_t^2\phi+s\mathcal O_{\lambda}(s\Delta x).
\end{align*}
Similarly, we have
\begin{align*}
{\bf A}_x q^{(11)}=\ & s\lambda\varphi\partial_t\phi+s\mathcal O_{\lambda}(s\Delta x),\\
\quad {\bf t}^+\big({\bf A}_x q^{(11)}\big)=\ & s\lambda\varphi\partial_t\phi+s\mathcal O_{\lambda}(\Delta t)+s\mathcal O_{\lambda}(s\Delta x),\\
{\bf D}_t{\bf A}_x q^{(11)} = \ &s\lambda^2{\bf A}_x\varphi|\partial_t\phi|^2+s\lambda{\bf A}_x\varphi\partial_t^2\phi+s\mathcal O_{\lambda}(s\Delta x)\\
=\ & s\lambda^2\varphi|\partial_t\phi|^2+s\lambda\varphi\partial_t^2\phi+s\mathcal O_{\lambda}((\Delta x)^2)+s\mathcal O_{\lambda}(s\Delta x),\\
{\bf t}^+\big({\bf D}_t{\bf A}_x q^{(11)}\big)=\ & s\lambda^2\varphi|\partial_t\phi|^2+s\lambda\varphi\partial_t^2\phi+s\mathcal O_{\lambda}(\Delta t)+s\mathcal O_{\lambda}((\Delta x)^2)+s\mathcal O_{\lambda}(s\Delta x).
\end{align*}
{Together with
\begin{align}\label{B.38}
\ &\mathbb E\int\limits_{\mathcal M^*} |{\bf D}_xY|^2\bigg|_{t=T}\nonumber\\
\displaybreak[0]\leq \  & s\mathcal O_{\lambda}(1)e^{s\mathcal O_{\lambda}(1)}\bigg(\mathbb E\int\limits_{\mathcal M^*} |{\bf D}_xy|^2\bigg|_{t=T}+\mathbb E\int\limits_{\mathcal M} |y|^2\bigg|_{t=T}\bigg)^{\frac{1}{2}}\bigg(\mathbb E\int\limits_{\mathcal M^*} |{\bf D}_xY|^2\bigg|_{t=T}\bigg)^{\frac{1}{2}}\nonumber\\
\leq \ & s\mathcal O_{\lambda}(1)e^{s\mathcal O_{\lambda}(1)}\bigg(\mathbb E\int\limits_{\mathcal M^*} |{\bf D}_xY|^2\bigg|_{t=T}\bigg)^{\frac{1}{2}},
\end{align}
}
due to the discrete energy estimate in Lemma A.1, we obtain the estimate for $D_9$ that
\begin{align}\label{}
		D_9		\ge\ &   \mathbb{E} \iint\limits_{\mathcal M^*\times\mathcal N}\left(s\lambda^2\varphi|\partial_t\phi|^2+s\lambda\varphi \partial_t^2\phi+s\mathcal O_{\lambda}(s\Delta x)\right)|{\bf D}_xY|^2+\nonumber\\
\ & \mathbb{E} \int\limits_{\mathcal M^*}\left(s\lambda \varphi \partial_t\phi-s\mathcal O_{\lambda}(s\Delta x)\right) |{\bf D}_xY|^2\bigg|_{t=0}-{s^2\mathcal O_{\lambda}(1)e^{s\mathcal O_{\lambda}(1)}\bigg(\mathbb E\int\limits_{\mathcal M^*} |{\bf D}_xY|^2\bigg|_{t=T}\bigg)^{\frac{1}{2}}}-\nonumber\\
\ & \mathcal O(\Delta t)\mathbb E\iint\limits_{\mathcal M^*\times\mathcal N^*}s\lambda\varphi|{\bf D}_xZ|^2.
	\end{align}

By an analogous argument to that used in estimating $q^{(11)}$, it follows that
\begin{align*}
q^{(12)}=-s\lambda\varphi\partial_x\phi+s\mathcal O_{\lambda}(s\Delta x),\quad {\bf D}_x q^{(12)}=-s\lambda^2\varphi|\partial_x\phi|^2-s\lambda\varphi\partial_x^2\phi+s\mathcal O_{\lambda}(s\Delta x).
\end{align*}
Then we have
\begin{align}
B_4\geq \ &\mathbb E\int\limits_{\mathcal N}\left(s\lambda\varphi\partial_x\phi+s\mathcal O_{\lambda}(s\Delta x)\right) {\bf tr}(|{\bf D}_xY|^2)\bigg |_{x=0}-\nonumber\\
\ &\mathbb E\int\limits_{\mathcal N}\left(s\lambda\varphi\partial_x\phi+s\mathcal O_{\lambda}(s\Delta x)\right) {\bf tr}(|{\bf D}_xY|^2)\bigg |_{x=1}
\end{align}
and
\begin{align}
D_{10}=\mathbb E\iint\limits_{\mathcal M^*\times \mathcal N}\left(s\lambda^2\varphi|\partial_x \phi|^2+s\lambda\varphi\partial_x^2\phi-s\mathcal O_{\lambda}(s\Delta x)\right)|{\bf D}_x Y|^2.
\end{align}

We now transfer to handle $D_{11}$ and $D_{12}$. Combining (\ref{A.4}) and Lemma 2.5,  we obtain
\begin{align*}
q^{(13)}=\ &s\lambda^2\varphi\left(|\partial_t\phi|^2-|\partial_x\phi|^2\right)+s\mathcal O_{\lambda}(s\Delta x)+(\Delta x)^2s\mathcal O_{\lambda}(1)s^2\mathcal O_{\lambda}(1)\nonumber\\
=\ &s\lambda^2\varphi\left(|\partial_t\phi|^2-|\partial_x\phi|^2\right)+s\mathcal O_{\lambda}(s\Delta x)=s\mathcal O_{\lambda}(1).
\end{align*}
and then
\begin{align*}
{\bf A}_x q^{(13)}=\ & s\lambda^2{\bf A}_x\left(\varphi\left(|\partial_t\phi|^2-|\partial_x\phi|^2\right)\right)+s\mathcal O_{\lambda}(s\Delta x)\nonumber\\
=\ & s\lambda^2\varphi\left(|\partial_t\phi|^2-|\partial_x\phi|^2\right)+s\mathcal O_{\lambda}((\Delta x)^2)+s\mathcal O_{\lambda}(s\Delta x)=s\mathcal O_{\lambda}(1).
\end{align*}
Further, by Lemma 2.5  and Lemma 2.6 we also have 
\begin{align*}
{\bf D}_x^2 q^{(13)}=\ & {\bf D}_x^2q^{(13)}_{1}{\bf A}_x^2q^{(13)}_{2}+2{\bf A}_x{\bf D}_xq^{(13)}_{1}{\bf A}_x{\bf D}_xq^{(13)}_{2}+{\bf A}_x^2q^{(13)}_{1}{\bf D}_x^2q^{(13)}_{2}=s\mathcal O_{\lambda}(1),
\end{align*}
where
\begin{align*}
q^{(13)}_{1}=\ &-r\varphi^{-1}\left({\bf D}_t\varphi{\bf D}_t\rho-{\bf D}_x\varphi {\bf D}_x\rho\right),\\
q^{(13)}_{2}=\ & 1+\frac{1}{2}(\Delta x)^2 r{\bf D}_x^2\rho 
\end{align*}
satisfy
\begin{align*}
{\bf A}_x^m{\bf D}_x^n q^{(13)}_{1}=\ &-\mathcal O_{\lambda}(1) \sum_{i+j=n}{\bf A}_x^i{\bf D}_x^j (r{\bf D}_t\rho)+\mathcal O_{\lambda}(1)\sum_{i+j=n} {\bf A}_x^i{\bf D}_x^j (r{\bf D}_x\rho)=s\mathcal O_{\lambda}(1),\\
{\bf A}_x^m{\bf D}_x^n q^{(13)}_{2}=\ &\mathcal O(1)+(\Delta x)^2{\bf A}_x^m{\bf D}_x^n(r{\bf D}_x^2\rho)=\mathcal O_{\lambda}(1).
\end{align*}
for $m,n\in \mathbb{N}$.
 Then we obtain 
\begin{align}
D_{11}\geq\ & -s\mathcal O_{\lambda}(1)\mathbb E\iint\limits_{\mathcal M\times \mathcal N}|Y|^2,\\
D_{12}\geq \ &-\mathbb{E} \iint\limits_{\mathcal M^*\times \mathcal N} \left(s\lambda^2\varphi\left(|\partial_t\phi|^2-|\partial_x\phi|^2\right)+s\mathcal O_{\lambda}(s\Delta x)\right)|{\bf D}_xY|^2.
\end{align}

Finally, combining the preceding estimates for $B_4$, $D_8, \cdots, D_{12}$ and noticing that $s\Delta x\leq \varepsilon$ and $\Delta t=\varepsilon\mathcal O((\Delta x)^2)$, we obtain the desired estimates (\ref{A.27}) and (\ref{A.28}) then complete the proof of Lemma B.3. \hfill$\Box$

\vspace{2mm}
	
{\noindent}{\bf Lemma B.4.}\ {\em Provided $s\Delta x\le \varepsilon$ and $\Delta t=\varepsilon\mathcal O\left((\Delta x)^2\right)$ for a sufficiently small $\varepsilon>0$ and $\beta\in (0,1)$, there exists positive constant $\lambda_0$ such that 
	\begin{align}\label{}
B_5\geq \ & \mathbb E\int\limits_{\mathcal M}\left(s^3\lambda^3\varphi^3\partial_t\phi(|\partial_t\phi|^2-|\partial_x\phi|^2)\right)|Y|^2\bigg |_{t=0}-{BT_4}
\end{align}
and
\begin{align}\label{A.40}
\sum_{j=13}^{17} D_j \ge\ & c_0\mathbb E\iint\limits_{\mathcal M\times\mathcal N}s^3\lambda^3\varphi^3|Y|^2-DT_3
\end{align}
for all fixed $\lambda\geq \lambda_0$, where $c_0=(4-4\beta)\min_{x\in I}|\partial_x\phi|^2$ and
\begin{align*}
DT_3=\ &  \left(\varepsilon s+1\right)s^2\mathcal O_{\lambda}(1)\mathbb E\iint\limits_{\mathcal M\times \mathcal N}|Y|^2+{\varepsilon s\mathcal O_{\lambda}(1)\mathbb E\iint\limits_{\mathcal M\times \mathcal N^*}|Z|^2}+\varepsilon s\mathcal O_{\lambda}(1)\mathbb E\iint\limits_{\mathcal M^*\times \mathcal N}|{\bf D}_xY|^2
\end{align*}
and

\begin{align*}
BT_4 =\ & \left(\varepsilon s+1\right)s^2\mathcal O_{\lambda}(1)\mathbb E\int\limits_{\mathcal M}|Y|^2\bigg |_{t=0}+ { s^3\mathcal O_{\lambda}(1)e^{s\mathcal O_{\lambda}(1)} \bigg(\mathbb E\int\limits_{\mathcal M}|Y|^2\bigg |_{t=T}\bigg)^{\frac{1}{2}}}.
\end{align*}
}

\vspace{2mm}

{\noindent\bf Proof.}\ For $q^{(21)}$, from Lemma 2.5 and
\begin{align}\label{A.41}
r{\bf D}_t^2 \rho-r{\bf D}_x^2\rho =\ &r\partial_t^2\rho+s^2\mathcal O_{\lambda}((s\Delta t)^2))-r\partial^2_x\rho+s^2\mathcal O_{\lambda}((s\Delta x)^2)\nonumber\\
=\ & (r\partial_t\rho)^2-(r\partial_x\rho)^2+s\mathcal O_{\lambda}(1)+s^2\mathcal O_{\lambda}(s\Delta x),
\end{align}
it follows that
\begin{align*}
q^{(21)}=\ & \left(r\partial_t\rho+s\mathcal O_{\lambda}((s\Delta t)^2)\right)\left((r\partial_t\rho)^2-(r\partial_x\rho)^2+s\mathcal O_{\lambda}(1)+s^2\mathcal O_{\lambda}(s\Delta x)\right)\\
=\ & -s^3\lambda^3\varphi^3\partial_t\phi(|\partial_t\phi|^2-|\partial_x\phi|^2)+s^2\mathcal O_{\lambda}(1)+s^3\mathcal O_{\lambda}(s\Delta x).
\end{align*} 
Further we have
\begin{align*}
{\bf t}^+\big(q^{(21)}\big)=\ & -s^3\lambda^3\varphi^3{\bf t}^+\left(\partial_t\phi(|\partial_t\phi|^2-|\partial_x\phi|^2)\right)+s^2\mathcal O_{\lambda}(1)+s^3\mathcal O_{\lambda}(s\Delta x)\\
=\ &-s^3\lambda^3\varphi^3\partial_t\phi(|\partial_t\phi|^2-|\partial_x\phi|^2)+s^3\mathcal O_{\lambda}(\Delta t)+s^2\mathcal O_{\lambda}(1)+s^3\mathcal O_{\lambda}(s\Delta x).
\end{align*}
Obviously, the same expression also holds for ${\bf t}^+\left({\bf t}^+(q^{21})\right)$. Similarly,
\begin{align*}
\ & {\bf D}_tq^{(21)}\\
=\ & {\bf D}_t\left( r{\bf A}_t{\bf D}_t\rho\right){\bf t}^+\left(r{\bf D}_t^2 \rho-r{\bf D}_x^2\rho\right)+{\bf t}^-\left(r{\bf A}_t{\bf D}_t\rho\right) {\bf D}_t\left(r{\bf D}_t^2 \rho-r{\bf D}_x^2\rho\right)\nonumber\\
=\ & \left(\partial_t(r\partial_t\rho)+s\mathcal O_{\lambda}((s\Delta t)^2)\right)\left((r\partial_t\rho)^2-(r\partial_x\rho)^2+s^2\mathcal O_{\lambda}(\Delta t)+s\mathcal O_{\lambda}(1)+s^2\mathcal O_{\lambda}(s\Delta x)\right)+\nonumber\\
\ &  \left(r\partial_t\rho+s\mathcal O_{\lambda}(\Delta t)\right)\left(\partial_t(r\partial_t^2\rho)-\partial_t(r\partial_x^2\rho)+s^2\mathcal O_{\lambda}((s\Delta t)^2)+s^2\mathcal O_{\lambda}(s\Delta x)\right)\nonumber\\
=\ & -3s^3\lambda^4\varphi^3|\partial_t\phi|^2(|\partial_t\phi|^2-|\partial_x\phi|^2)-s^3\lambda^3\varphi^3\partial_t^2\phi
(3|\partial_t\phi|^2-|\partial_x\phi|^2)+s^2\mathcal O_{\lambda}(1)+s^3\mathcal O_{\lambda}(s\Delta x).\end{align*}
The we obtain
\begin{align}\label{B.47}
 B_5 \geq \ & \mathbb E\int\limits_{\mathcal M}\left(s^3\lambda^3\varphi^3\partial_t\phi(|\partial_t\phi|^2-|\partial_x\phi|^2)-s^2\mathcal O_{\lambda}(1)-s^3\mathcal O_{\lambda}(\Delta x)\right)|Y|^2\bigg |_{t=0}-\nonumber\\
\ & s^3\mathcal O_{\lambda}(1) \mathbb E\int\limits_{\mathcal M}|Y|^2\bigg |_{t=T}\nonumber\\
\geq \  & \mathbb E\int\limits_{\mathcal M}\left(s^3\lambda^3\varphi^3\partial_t\phi(|\partial_t\phi|^2-|\partial_x\phi|^2)-s^2\mathcal O_{\lambda}(1)-s^3\mathcal O_{\lambda}(\Delta x)\right)|Y|^2\bigg |_{t=0}-\nonumber\\
\ &{ s^3\mathcal O_{\lambda}(1)e^{s\mathcal O_{\lambda}(1)} \bigg(\mathbb E\int\limits_{\mathcal M}|Y|^2\bigg |_{t=T}\bigg)^{\frac{1}{2}}}
\end{align}
and
\begin{align}
\label{A.43}D_{13}\geq\ & \mathbb E\iint\limits_{\mathcal M\times\mathcal N}\Big(3s^3\lambda^4\varphi^3|\partial_t\phi|^2(|\partial_t\phi|^2-|\partial_x\phi|^2)+s^3\lambda^3\varphi^3\partial_t^2\phi(3|\partial_t\phi|^2-|\partial_x\phi|^2)
-\nonumber\\
\ &s^3\mathcal O_{\lambda}(\Delta t)-s^2\mathcal O_{\lambda}(1)-s^3\mathcal O_{\lambda}(s\Delta x)\Big)|Y|^2,\\
D_{14}\geq \ & {-s\mathcal O_{\lambda}(s^2\Delta t)\mathbb E\iint\limits_{\mathcal M\times \mathcal N^*}|Z|^2,}
\end{align}
{Here the last inequality in (\ref{B.47}) follows from 
\begin{align*}
\ &\mathbb E\int\limits_{\mathcal M} |Y|^2\bigg|_{t=T}\leq e^{s\mathcal O_{\lambda}(1)}\bigg(\mathbb E\int\limits_{\mathcal M} |Y|^2\bigg|_{t=T}\bigg)^{\frac{1}{2}},
\end{align*}
which is analogous to (\ref{B.38}).
}

Using (\ref{A.41}),  Lemma 2.5 and Lemma 2.6, we obtain the expression for ${\bf D}_x q^{(22)}$ that
\begin{align*}
{\bf D}_xq^{(22)}= \ &  -{\bf D}_x\left( r{\bf A}_x{\bf D}_x\rho\right){\bf A}_x\left(r{\bf D}_t^2 \rho-r{\bf D}_x^2\rho\right)-{\bf A}_x\left(r{\bf A}_x{\bf D}_x\rho\right) {\bf D}_x\left(r{\bf D}_t^2 \rho-r{\bf D}_x^2\rho\right)\\
\displaybreak[0]=\ & -\left(\partial_x(r\partial_x\rho)+s\mathcal O_{\lambda}((s\Delta x)^2)\right)\left((r\partial_t\rho)^2-(r\partial_x\rho)^2+s\mathcal O_{\lambda}(1)+s^2\mathcal O_{\lambda}(s\Delta x)\right)-\\
\ &  \left(r\partial_x\rho+s\mathcal O_{\lambda}((s\Delta x)^2)\right)\left(\partial_x(r\partial_t^2\rho)-\partial_x(r\partial_x^2\rho)+s^2\mathcal O_{\lambda}(s\Delta x)\right)\\
=\ & 3s^3\lambda^4\varphi^3|\partial_x\phi|^2(|\partial_t\phi|^2-|\partial_x\phi|^2)+s^3\lambda^3\varphi^3\partial_x^2\phi
(|\partial_t\phi|^2-3|\partial_x\phi|^2)
+\\
\ & s^2\mathcal O_{\lambda}(1)+s^3\mathcal O_{\lambda}(s\Delta x).\end{align*}
The expression of ${\bf A}_x{\bf D}_xq^{(22)}$ is the same as the one of ${\bf D}_xq^{(22)}$. Therefore, we obtain the following estimates for $D_{15}$ and $D_{16}$:
\begin{align}
D_{15}\geq\ & \mathbb E\iint\limits_{\mathcal M\times\mathcal N}\Big(-3s^3\lambda^4\varphi^3|\partial_x\phi|^2(|\partial_t\phi|^2-|\partial_x\phi|^2)-s^3\lambda^3\varphi^3\partial_x^2\phi
(|\partial_t\phi|^2-3|\partial_x\phi|^2)-\nonumber\\
\ &s^2\mathcal O_{\lambda}(1)-s^3\mathcal O_{\lambda}(s\Delta x)\Big)|Y|^2,\\
D_{16}\geq\ & -s\mathcal O_{\lambda} \left((s\Delta x)^2\right)\mathbb E\iint\limits_{\mathcal M^*\times \mathcal N}|{\bf D}_x Y|^2. 
\end{align}

By (\ref{A.4}) and (\ref{A.41}), we easily see that
\begin{align*}
 q^{(23)}
= \ &\left(-s\lambda^2\varphi\left(|\partial_t\phi|^2-|\partial_x\phi|^2\right)\right)\left(s^2\lambda^2\varphi^2(|\partial_t\phi
|^2-|\partial_x\phi|^2)\right)+s^2\mathcal O_{\lambda}(1)+s^3\mathcal O_{\lambda}\left(s\Delta x\right)\\
=\ &-s^3\lambda^4\varphi^3\left(|\partial_t\phi|^2-|\partial_x\phi|^2\right)^2+s^2\mathcal O_{\lambda}(1)+s^3\mathcal O_{\lambda}\left(s\Delta x\right)
\end{align*}
and then
\begin{align}
\label{A.47}D_{17}\geq  -\mathbb{E}\iint\limits_{\mathcal M\times\mathcal N}\left(s^3\lambda^4\varphi^3\left(|\partial_t\phi|^2-|\partial_x\phi|^2\right)^2+s^2\mathcal O_{\lambda}(1)+s^3\mathcal O_{\lambda}\left(s\Delta x\right)\right)|Y|^2.
\end{align}

	Combining (\ref{A.43})-(\ref{A.47}), we obtain
\begin{align}\label{A.48}
\sum_{j=13}^{17} D_j\geq \ & \mathbb E\iint\limits_{\mathcal M\times\mathcal N}s^3\lambda^3\varphi^3\left(2\lambda H_{s,\lambda}^2(\phi)-(6\beta+2)H_{s,\lambda}(\phi)+(4-4\beta)|\partial_x\phi|^2\right)|Y|^2-\nonumber\\
\ & \left(s^2\mathcal O_{\lambda}(1)+s^3\mathcal O_{\lambda}(s\Delta x)\right)\mathbb E\iint\limits_{\mathcal M\times \mathcal N}|Y|^2-s\mathcal O_{\lambda}(s^2\Delta t)\mathbb E\iint\limits_{\mathcal M\times \mathcal N}{\bf t}^+(|Z|^2)-\nonumber\\
\ & s\mathcal O_{\lambda}(s\Delta x)\mathbb E\iint\limits_{\mathcal M^*\times \mathcal N}|{\bf D}_xY|^2
\end{align}
with
\begin{align*}
H_{s,\lambda}(\phi)=|\partial_t\phi|^2-|\partial_x\phi|^2.
\end{align*}
Then we can choose $\lambda\geq \lambda_0$ with sufficiently $\lambda_0$ to obtain 
\begin{align}\label{A.49}
2\lambda H_{s,\lambda}^2(\phi)-(6\beta+2)H_{s,\lambda}(\phi)+(4-4\beta)|\partial_x\phi|^2>0 
\end{align}
for $H_{s,\lambda}(\phi)\neq 0$. If $H_{s,\lambda}(\phi)= 0$, this inequality (\ref{A.49}) also holds due to $\beta\in (0,1)$. Consequently, there  exist positive constant $\lambda_0$ such that for fixed $\lambda>\lambda_0$, we have
\begin{align}\label{A.50}
\ &\mathbb E\iint\limits_{\mathcal M\times\mathcal N}s^3\lambda^3\varphi^3\left(2\lambda H_{s,\lambda}^2(\varphi)-(6\beta+2)H_{s,\lambda}(\varphi)+(4-4\beta)|\partial_x\phi|^2\right)|Y|^2\nonumber\\
\geq\ & c_0\mathbb E\iint\limits_{\mathcal M\times\mathcal N}s^3\lambda^3\varphi^3|Y|^2.
\end{align} 
Finally, substituting (\ref{A.50}) into (\ref{A.48}), we obtain the desired estimate (\ref{A.40}). The proof of Lemma B.4 is completed.\hfill$\Box$

\vspace{5mm}

{\noindent\bf Acknowledgement}

\vspace{2mm}

\noindent The first author is supported by National Natural Science Foundation of China (No. 12171248) and The Natural Science Foundation of Jiangsu Province (No. BK20230412).

\vspace{2mm}

{\noindent\bf Data availability statement}

\vspace{2mm}

\noindent All data that support the findings of this study are included within the article (and any
supplementary files).

\vspace{5mm}

	\newcounter{cankao}
	\begin{list}
		{[\arabic{cankao}]}{\usecounter{cankao}\itemsep=0cm} \centerline{\bf
			References} \vspace*{0.5cm} \small


\item\label{Bao-SIAMNA-2020} Bao G, Lin Y and Xu X. Inverse scattering by a random periodic structure. SIAM Journal on Numerical Analysis, 2020, 58(5): 2934-2952.

\item\label{Bao2012} Bao G and Xu X. An inverse random source problem in quantifying the elastic modulus of nanomaterials. Inverse Problems, 2012, 29(1): 015006.

\item\label{Baudouin-SICON-2013}Baudouin L and Ervedoza S. Convergence of an inverse problem for a 1-D discrete wave equation. SIAM Journal on Control and Optimization, 2013, 51(1): 556-598.

    \item\label{Bau2015} Baudouin L, Ervedoza S and Osses A. Stability of an inverse problem for the discrete wave equation and convergence results. Journal de Mathématiques Pures et Appliquées, 2015, 103(6): 1475-1522.
    
\item\label{Bellassoued-Springer-2017}    Bellassoued M and Yamamoto M. Carleman estimates and applications to inverse problems for hyperbolic systems. Tokyo: Springer, 2017.

\item\label{Boyer-ESAIM:COCV-2005} Boyer F and Hern\'{a}ndez-Santamar\'{i}a V. Carleman estimates for time-discrete parabolic equations and applications to controllability.   ESAIM: Control,  Optimisation and Calculus of Variations, 2020, 26: 12.
    	
\item\label{Boyer-JMPA-2010}  Boyer F, Hubert F and Le Rousseau J. Discrete carleman estimates for elliptic operators and uniform controllability of semi-discretized parabolic equations.  Journal de math\'{e}matiques pures et appliqu\'{e}es, 2010, 93(3): 240-276.

\item\label{Boyer-POINCARE-AN-2014} Boyer F and Le Rousseau J. Carleman estimates for semi-discrete parabolic operators and application to the controllability of semi-linear semi-discrete parabolic equations. Annales de l'Institut Henri Poincar\'{e}, Analyse non lin\'{e}aire, 2014, 31(5): 1035-1078.
    
{    
 \item\label{Bukhgeim-SM-1981} Bukhgeim A L and  Klibanov M V, Uniqueness in the large of a class of multidimensional inverse problems, Soviet Mathematics-Doklady, 1981, 24: 244-247.
 }
    
\item\label{Carreno-ACM-2023}   Casanova P G and Hern\'{a}ndez-Santamar\'{i}a V. Carleman estimates and controllability results for fully discrete approximations of 1D parabolic equations. Advances in Computational Mathematics, 2021, 47:  1-71.
    
\item\label{Cav 2002} Cavalier L and Tsybakov A. Sharp adaptation for inverse problems with random noise. Probability Theory and Related Fields, 2002, 123(3): 323-354

\item\label{Cerpa-JMPA-2022}  Cerpa E, Lecaros R, Nguyen T and P\'{e}rez A. Carleman estimates and controllability for a semi-discrete fourth-order parabolic equation. Journal de Math\'{e}matiques Pures et Appliqu\'{e}es, 2022, 164: 93-130.
    
{    
\item\label{Clark-DAM-1987}    Clark D S. Short proof of a discrete Gronwall inequality. Discrete Applied Mathematics, 1987, 16(3): 279-281.
}
    
\item\label{Dou-IP-2024}     Dou F, Lü P and Wang Y. Stability and regularization for ill-posed Cauchy problem of a stochastic parabolic differential equation. Inverse Problems, 2024, 40(11): 115005.
    
\item\label{Hernández2021} Hernández-Santamaría V. Controllability of a simplified time-discrete stabilized Kuramoto-Sivashinsky system. Evolution Equations and Control Theory, 2023, 12(2): 459-501.
    
\item\label{Imanuvilov-IP-2021} Imanuvilov O Y and Yamamoto M. GlobalLipschitz stability in an inverse hyperbolic problem by interior observations. Inverse problems, 2001, 17(4): 717.
    
 \item\label{Klibanov-IP-2024}   Klibanov M V. A new type of ill-posed and inverse problems for parabolic equations. Communications on Analysis and Computation, 2024, 2(4): 367-398.
{     
\item\label{Klibanov-DG-2021}  Klibanov M V and  Li J. Inverse problems and Carleman estimates: global
uniqueness, global convergence and experimental data. De Gruyter, 2021.

\item\label{Klibanov-IPI-2021}Klibanov M V, Li J, Yang Z. Convexification numerical method for a coefficient inverse problem for the system of nonlinear parabolic equations governing mean field games. Inverse Problems and Imaging, 2025, 19(2): 219-252.
}

\item\label{Kumar-2024}   Kumar M. Boundary controllability of a fully-discrete approximation of linear KdV equation. preprint, 2025. https://hal.science/hal-04828698v2
    
{
 \item\label{Lecaros-IP-2025}    Lecaros R, López-Ríos J and Pérez A A. Lipschitz stability and reconstruction in inverse problems for semi-discrete parabolic operators. Inverse Problems, 2025, 41(11): 115012.
}
     
 \item\label{Lecaros-JDE-2023}  Lecaros R, Morales R, P\'{e}rez A and Zamorano S. Discrete carleman estimates and application to controllability for a fully-discrete parabolic operator with dynamic boundary conditions. Journal of Differential Equations, 2023, 365:  832-881.
     
\item\label{Lecaros-ESIAM-2021}  Lecaros R, Ortega J H and Pérez A. Stability estimate for the semi-discrete linearized Benjamin-Bona-Mahony equation. ESAIM: Control, Optimisation and Calculus of Variations, 2021, 27: 93.

\item\label{Li-Xu-SIAMUQ-2022}     Li P and Wang X. An inverse random source problem for the biharmonic wave equation. SIAM/ASA Journal on Uncertainty Quantification, 2022, 10(3): 949-974.

 \item\label{Lu-IP-2012}    Lü Q. Carleman estimate for stochastic parabolic equations and inverse stochastic parabolic problems. Inverse Problems, 2012, 28(4): 045008.
     
 \item\label{Lu-IP-2013} Lü Q. Observability estimate and state observation problems for stochastic hyperbolic equations. Inverse Problems, 2013, 29(9): 095011.

    \item\label{Lu-Wang-2025} Lü Q and Wang Y. An inverse source problem for semilinear stochastic hyperbolic equations. arXiv preprint arXiv:2504.17398, 2025.
	
\item\label{Lü 2015} Lü Q and Zhang X. Global uniqueness for an inverse stochastic hyperbolic problem with three unknowns. Communications on Pure and Applied Mathematics, 2015, 68(6): 948-63.	
    
    \item\label{Lu2023} Lü Q and Zhang X. Inverse problems for stochastic partial differential equations: some progresses and open problems. Numerical Algebra, Control and Optimization, 2024, 14(2): 227-272.
        
 \item\label{Wang-arXiv-2024}  Wang Y and Zhao Q. Null controllability for stochastic fourth order semi-discrete parabolic equations. arXiv preprint arXiv:2405.03257, 2024.
    
\item\label{Wu 2022} Wu B and Liu J. On the stability of recovering two sources and initial status in a stochastic hyperbolic-parabolic system. Inverse Problems, 2022, 38(2): 025010.
	
	\item\label{Wu 2024} Wu B, Wang Y and Wang Z. Carleman estimates for space semi-discrete approximations of one-dimensional stochastic parabolic equation and its applications. Inverse Problems, 2024, 40(11): 115003.	

{
\item\label{Yamamoto-SP-2025} Yamamoto M. Introduction to Inverse Problems for Evolution Equations: Stability and Uniqueness by Carleman Estimates. Springer Nature, 2025.
}

\item\label{Yu-MMAS-2023} Yu Y and Zhao Q. Global uniqueness in an inverse problem for a class of damped stochastic plate equations. Mathematical Methods in the Applied Sciences, 2023, 46(1): 765-776.

\item\label{Yuan-IP-2017} Yuan G. Conditional stability in determination of initial data for stochastic parabolic equations. Inverse Problems, 2017, 33(3): 035014.

\item\label{Yuan-JMAA-2017} Yuan G. Determination of two unknowns simultaneously for stochastic Euler–Bernoulli beam equations. Journal of Mathematical Analysis and Applications, 2017, 450(1): 137-151.    
    
\item\label{Zhang 2022} Zhang W and Zhao Z. Convergence analysis of a coefficient inverse problem for the semi-discrete damped wave equation. Applicable Analysis, 2022, 101(4): 1430-1455.

\item\label{Zha-SIMA-2008}Zhang X. Carleman and observability estimates for stochastic wave equations. SIAM Journal on Mathematical Analysis, 2008, 40(2): 851-868.
    
\item\label{Zhao-SICON-2025} Zhao Q. Null Controllability for Stochastic Semidiscrete Parabolic Equations. SIAM Journal on Control and Optimization, 2025, 63(3): 2007-2028.

\item\label{Zuazua-SIMA-2005}Zuazua E. Propagation, observation, and control of waves approximated by finite difference methods. SIAM review, 2005, 47(2): 197-243.

	\end{list}
\end{document}